\documentclass[12pt]{article}
\usepackage{theorem}
\usepackage{amsmath,amssymb,amscd,latexsym,eufrak}

\newcommand{\bvskip}{\vspace{7mm}}
\newcommand{\mvskip}{\vspace{5mm}}
\newcommand{\svskip}{\vspace{3mm}}
\newcommand{\bhskip}{\hspace{1cm}}
\newcommand{\mhskip}{\hspace{5mm}}
\newcommand{\shskip}{\hspace{22pt}}
\newcommand{\no}[1]{\svskip \noindent {\bf {#1}.}\ }
\newcommand{\C}{{\Bbb C}}
\newcommand{\R}{{\Bbb R}}
\newcommand{\N}{{\Bbb N}}
\newcommand{\Z}{{\Bbb Z}}
\newcommand{\BP}{{\Bbb P}}
\newcommand{\Q}{{\Bbb Q}}
\newcommand{\Supp}{{\rm Supp}\:}
\newcommand{\Sing}{{\rm Sing}\:}
\newcommand{\Bs}{{\rm Bs}\:}
\newcommand{\Bk}{{\rm Bk}\:}
\newcommand{\res}{{\rm res}\:}
\newcommand{\red}{{\rm red}\:}
\newcommand{\Tang}{{\rm Tang}\,}
\newcommand{\QED}{{\unskip\nobreak\hfil\penalty50\quad\null\nobreak\hfil
{ Q.E.D.}\parfillskip0pt\finalhyphendemerits0\par\medskip}}
\newcommand{\Proof}{\noindent{\bf Proof.}\quad}
\newcommand{\A}{{\bf A}}
\newcommand{\Spec}{{\rm Spec}\:}
\newcommand{\Proj}{{\rm Proj}\:}
\newcommand{\lto}{\longrightarrow}
\newcommand{\dlto}{\:\cdots\!\!\to}
\newcommand{\correspond}[1]{\raisebox{-1.7ex}{$\,\,\buildrel\displaystyle
{\sim}\over{\scriptstyle{#1}}\,\,$}}
\newcommand{\Kdim}{\mbox{{\rm K-dim}}\:}
\newcommand{\lkd}{\ol{\kappa}}
\newcommand{\kd}{\kappa}
\newcommand{\ord}{{\rm ord}}
\newcommand{\codim}{{\rm codim}\:}
\newcommand{\height}{{\rm ht}\:}
\newcommand{\length}{{\rm length}\:}
\newcommand{\trdeg}{{\rm tr.deg}\:}
\newcommand{\trace}{{\rm Tr}\:}
\newcommand{\rank}{{\rm rank}\:}
\newcommand{\Aut}{{\rm Aut}\:}
\newcommand{\Pic}{{\rm Pic}\:}
\newcommand{\Alb}{{\rm Alb}\:}
\newcommand{\Cl}{{\rm C\ell}\:}
\newcommand{\GL}{{\rm GL}\:}
\newcommand{\SSL}{{\rm SL}\:}
\newcommand{\PGL}{{\rm PGL}\:}
\newcommand{\Gal}{{\rm Gal}\:}
\newcommand{\Div}{{\rm Div}\:}
\newcommand{\divisor}{{\rm div}\:}
\newcommand{\diag}{{\rm diag}\:}
\newcommand{\NS}{{\rm NS}\:}
\newcommand{\NE}{{\rm NE}\:}
\newcommand{\Der}{{\rm Der}\:}
\newcommand{\indlim}[1]{\displaystyle{\lim_
{\textstyle{\lto}\atop{\scriptstyle{#1}}}}}
\newcommand{\Ker}{{\rm Ker}\:}
\newcommand{\Coker}{{\rm Coker}\:}
\newcommand{\im}{{\rm Im}\:}
\newcommand{\id}{{\rm id}\:}
\newcommand{\gr}{{\rm gr}\,}
\newcommand{\propersubset}{\raisebox{-.7ex}{$\,\,\buildrel\displaystyle
{\subset}\over{\scriptstyle\neq}\,\,$}}
\newcommand{\propersupset}{\raisebox{-.7ex}{$\,\,\buildrel\displaystyle
{\supset}\over{\scriptstyle\neq}\,\,$}}
\newcommand{\Cech}{$\stackrel{\scriptscriptstyle{\vee}}{\rm C}$ech~}
\newcommand{\SP}{{\mathcal P}}
\newcommand{\SQ}{{\mathcal Q}}
\newcommand{\SR}{{\mathcal R}}
\newcommand{\CS}{{\mathcal S}}
\newcommand{\ST}{{\mathcal T}}
\newcommand{\SL}{{\mathcal L}}
\newcommand{\SM}{{\mathcal M}}
\newcommand{\SN}{{\mathcal N}}
\newcommand{\SE}{{\mathcal E}}
\newcommand{\SF}{{\mathcal F}}
\newcommand{\SG}{{\mathcal G}}
\newcommand{\SH}{{\mathcal H}}
\newcommand{\SI}{{\mathcal I}}
\newcommand{\SJ}{{\mathcal J}}
\newcommand{\SK}{{\mathcal K}}
\newcommand{\SC}{{\mathcal C}}
\newcommand{\SA}{{\mathcal A}}
\newcommand{\SB}{{\mathcal B}}
\newcommand{\SO}{{\mathcal O}}
\newcommand{\SU}{{\mathcal U}}
\newcommand{\SV}{{\mathcal V}}
\newcommand{\SW}{{\mathcal W}}
\newcommand{\ga}{\EuFrak{a}}
\newcommand{\gb}{\EuFrak{b}}
\newcommand{\gm}{\EuFrak{m}}
\newcommand{\gn}{\EuFrak{n}}
\newcommand{\gp}{\EuFrak{p}}
\newcommand{\gq}{\EuFrak{q}}
\newcommand{\GA}{\EuFrak{A}}
\newcommand{\GB}{\EuFrak{B}}
\newcommand{\GE}{\EuFrak{E}}
\newcommand{\GM}{\EuFrak{M}}
\newcommand{\GN}{\EuFrak{N}}
\newcommand{\GP}{\EuFrak{P}}
\newcommand{\GQ}{\EuFrak{Q}}
\newcommand{\GS}{\EuFrak{Ye}}
\newcommand{\GU}{\EuFrak{U}}
\newcommand{\GV}{\EuFrak{V}}
\newcommand{\GW}{\EuFrak{W}}
\newcommand{\CH}[1]{\stackrel{\scriptscriptstyle{\vee}}{H}\!{}^{#1}}
\newcommand{\Hom}{{\rm Hom}}
\newcommand{\SHom}{{\cal H}{\em om}}
\newcommand{\Ext}{{\rm Ext}}
\newcommand{\st}[1]{\stackrel{{#1}}{\longrightarrow}}
\newcommand{\nwa}[1]{\nwarrow{\scriptstyle{#1}}}
\newcommand{\pbx}[2]{\parbox[t]{3.5cm}{#1}\ \hspace{5mm}\ \parbox[t]{7.8cm}{#2}
\vspace{0mm}}
\newcommand{\scs}[1]{\scriptstyle{#1}}
\newcommand{\dps}[1]{\displaystyle{#1}}
\newcommand{\mpg}[2]{\noindent\begin{minipage}[t]{#1}{#2}\end{minipage}\svskip}
\newcommand{\wt}{\widetilde}
\newcommand{\ol}{\overline}
\newcommand{\wh}{\widehat}
\newcommand{\quot}{/\!/}
\newcommand{\ch}{{\rm char\,}}
\newcommand{\is}[2]{({#1}\cdot{#2})}
\newcommand{\sis}[1]{({#1}^2)}
\newtheorem{thm}{Theorem}
\newtheorem{lem}[thm]{Lemma}

\begin{document}
\title{Equivariant classification of Gorenstein open log del Pezzo surfaces 
with finite group actions}
\author{Masayoshi Miyanishi and De-Qi Zhang}
\date{}
\maketitle

\begin{abstract}
We classify equivariantly Gorenstein log del Pezzo surfaces with boundaries at infinity and 
with finite group actions such that the quotient surface modulo the finite group has Picard 
number one. We also determine the corresponding finite groups. Better figures are available
upon request.
\end{abstract}

\setcounter{section}{-1}
\section{Introduction}
Let $\ol{X}$ be a normal projective rational surface with at worst rational double singularities 
and let $\ol{D}$ be a reduced divisor on $\ol{X}$. We, further, assume that there is 
a finite group $G$ acting faithfully on $\ol{X}$ so that $\ol{D}$ is $G$-stable.
We assume that $(\ol{X},\ol{D})$ has log terminal singularities 
(cf. \cite{KMM}, \cite{Ko}, \cite{oas}) and 
$\lkd(\ol{X} \setminus \ol{D}) = -\infty$. Let $f : X \to \ol{X}$ be the minimal resolution. 
Let $D$ be the proper transform of $\ol{D}$. We can write $f^*(\ol{D})=D+\Delta$, where $\Delta$ 
is a positive $\Q$-divisor such that $\Supp(\Delta)$ is the exceptional locus of $f$ 
arising only from the singular points lying on $\ol{D}$ (cf. {\em ibid.}). It is also known 
(cf. {\em ibid.}) that the exceptional graph of $f^{-1}(P)$ with $P\in \ol{D}\cap \Sing(\ol{X})$ 
is a linear $(-2)$-chain and that one of the end components of the chain meets transversally $D$ in 
one point. The $G$-action on $\ol{X}$ lifts up to a $G$-action on $X$ such that $D+\Delta$ 
is $G$-stable.

\begin{flushleft}
\underline{\hspace{10cm}}
\end{flushleft}

{\em 2000 Mathematics Subject Classification}\ \  Primary 14J17, 14J50; Secondary 14J26

{\em Key words and Phrases}\ \ Gorenstein singularity, open log del Pezzo surface, finite group 
action

\newpage
Our objective is to describe a pair $(\ol{X},\ol{D})$ with finite group action of $G$. 
In the present article we shall determine the geometric structure of 
the minimal resolution 
$(X, D+\Delta_\red + A)$ (see Remark (2) below)
of $(\ol{X},\ol{D})$ as well as the group action $G$ on ${\ol X}$.
We assume the following:

\svskip \noindent
{\bf Hypothesis (H)}\ \  {\em ${\ol X}$ has at worst rational double singularities,
$({\ol X}, {\ol D})$ is log terminal, ${\ol D} \ne 0$, $\kappa({\ol X} \setminus {\ol D}) =
- \infty$ and $\rho({\ol X} \quot G) = 1$.}

\svskip \noindent
{\bf Theorem A}\ \ {\em Assume the Hypothesis (H). Then either $K_X^2 \ge 8$
and one of the cases $(1a) - (1f)$ in Lemma $4$ occurs, or $2 \le K_X^2 \le 6$
and $f^{-1}{\ol D} + f^{-1}(\Sing {\ol X}) = D + \Delta_{\red} + A$
{\em (see Remark $(2)$ below)} is given in Figure m for some $1 \le m \le 43$ 
{\em (see Section 2)}.}

\svskip \noindent
{\bf Remark}\ \ (1) In \cite{Z2}, the equivariant classification of the pair 
$(X, G)$ where $X$ is smooth, is treated. In \cite{KM}, the authors dealt with the pair
$({\ol X}, G)$ where $\rho({\ol X}\quot G) \ge 2$, but
${\ol X}$ is assumed to be only log terminal. In \cite{Z3}, 
a finiteness criterion for $|Aut(X)|$ is given, where $X$ is Gorenstein
del Pezzo of Picard number one. See also \cite{KZ}.

\svskip \noindent
(2) We can write $f^{-1}({\ol D}) = D + \Delta_{\red}$
and $f^{-1}(\Sing {\ol X}) = \Delta_{\red} + A$, where
$A$ is contractible to singular points on ${\ol X}$ but not
on ${\ol D}$. We let $A = \sum_i A_i$ be the irreducible decomposition.

\svskip \noindent
(3) Note that either $D$ is irreducible or $D = D_1 + D_2$ is a linear
chain of two smooth rational curves (see Lemma 3). 
Figure $m$ contains the graph of $D + \Delta_{\red} + A$
$+$ (some $(-1)$-curves like $E, E_i, F_j$),
where each $(-1)$-curve (resp. each other curve) is represented by a 
broken (resp. solid line) with self intersection typed next to it.

\svskip \noindent
(4) Also shown in each Figure $m$ is a $\BP^1$-fibration 
$\Phi : X \rightarrow B$ ($\cong \BP^1$)
with all its singular fibres drawn vertically; thus one can read off
from each Figure $m$, the Picard number $\rho(X)$ and $K_X^2$ 
by blowing down $X$ to a Hirzebruch surface. See Lemma 14.

\svskip\noindent
(5) We use the notation $\Delta = \sum_{i=1}^t \Delta_i,
(\Delta_i)_{red} = \sum_{j=1}^{s_i} \Delta_i(j)$ in Lemma 2.
In each Figure $m$, if $\Delta$ or $\Delta_i$ is irreducible
we use the same letter to denote its support 
(which is a reduced irreducible curve).
If $t = 1$, we set $\Delta(j) := \Delta_1(j)$.

\svskip \noindent
{\bf Theorem B}\ \ {\em Concerning the group $G$ acting faithfully on $\ol{X}$ (or $X$), 
the following assertions hold:
\begin{enumerate}
\item[{\rm (1)}]
With the notations and assumptions in Theorem A, each
$G$ is determined in Section $2$.

\item[{\rm (2)}]
Conversely, given Figure $m$ of $D + \Delta_\red +A$
on a surface $X$ for some $1 \le m \le 43$ in Section 2,
we let $f : X \rightarrow {\ol X}$ be the contraction of $\Delta + A$
and ${\ol D} = f_*D$. Then we can find a finite group $G$
specified in $\S 2$ acting on ${\ol X}$ faithfully such that
the Hypothesis (H) is satisfied.
\end{enumerate}}

\svskip \noindent
{\bf Theorem C.} {\em With the notations and assumptions in Theorem A,
assume further that $K_X^2 \le 4$. Then either $G$ is soluble or
$|G:H| \le 2$ for some group $H$ in
$\{g = (a_{ij})\in \PGL_2({\C}) \mid a_{ij} \ne 0$ only when 
$i = j$ or $(i,j) = (3,1) \}$. We have $G = H$ except the case of
Figure 25.
\par

Conversely, any finite group in $PGL_2({\Bbb C})$ of the form above
can act on some ${\ol X}$ faithfully so that the Hypothesis (H)
is satisfied.}
\svskip

W assume throughout the article that the ground field $k$ is an algebraically closed 
field of characteristic zero.
\svskip

\noindent
{\bf Acknowledgement} \ \ 
The authors would like to thank the referee for the careful reading and
suggestions which improve the paper.
The first author would like to express his gratitude to 
National University of Singapore for the 
financial support under the program \lq\lq 
Scheme for Promotion of Research Interaction '' 
which realized this joint work. The second author was supported by
an Academic Research Fund of NUS. The authors like to thank the referee
for careful reading and suggestions which improve the paper.

\section{Geometric structure of the surface $X$}

Let us begin with the following result.
We assume Hypothesis (H) in \S 1.

\begin{lem}
The following conditions are equivalent:
\begin{enumerate}
\item[{\rm (1)}]
The Picard number $\rho(\ol{X}\quot G)=1$.
\item[{\rm (2)}]
$(\Pic\ol{X})^G\cong \Z$.
\item[{\rm (3)}]
$(\Pic\ol{X})^G\otimes\Q=((\Pic\ol{X})\otimes\Q)^G\cong \Q$.
\end{enumerate}
\end{lem}

\Proof
Since the pull back of the quotient map ${\ol X} \rightarrow {\ol X} \quot G$
induces an isomorphism between (Pic(${\ol X} \quot G)) \otimes {\Bbb Q}$ and
(Pic${\ol X})^G \otimes {\Bbb Q}$, (1) and (3) are equivalent.
Since ${\ol X}$ has at worst quotient singularities,
the resolution $f$ induces an isomorphism $\pi_1(X) \rightarrow \pi_1({\ol X})$
by Theorem 7.8 in \cite{Ko2}. So $\pi_1(\ol X) = 1$ and hence $\Pic {\ol X}$
and $(\Pic {\ol X})^G$ have no torsion elements; thus (2) and (3) are equivalent.
This proves the lemma.
\QED

From now on, we assume one of the equivalent conditions of Lemma 1. We call such 
a pair $(\ol{X},\ol{D})$ with a finite group action of $G$ a {\em  Gorenstein open 
log del Pezzo surface} provided $\ol{D} \ne 0$.

\begin{lem}
The following assertions hold.
\begin{enumerate}
\item[{\rm (1)}]
$\ol{D}, -K_{\ol{X}}$ and $-(K_{\ol{X}}+\ol{D})$ are all $\Q$-ample divisors. Each of 
them generates $(\Pic\ol{X})\otimes\Q$. The divisors
$D+\Delta, -K_X$ and $-(K_X+D+\Delta)$ are all nef and big and $\Q$-proportional 
to one another. 
\item[{\rm (2)}]
$f : X \rightarrow \ol{X}$ is nothing but the contraction of all $(-2)$-curves on $X$.
If $C$ is a curve on $X$ with $C^2 < 0$, then $C$ is either a $(-1)$-curve
or a $(-2)$-curve. If $X \rightarrow \Sigma_d$ is a birational morphism
to a Hirzebruch surface of degree $d$, then $d = 0, 1, 2$.

\item[{\rm (3)}] If $\Phi : X \rightarrow {\BP}^1$ is a $\BP^1$-fibration
and $\Gamma_1$ a singular fibre, then either (type $I_n$) 
$\Gamma_1 = E_1 + A_1 + \cdots + A_n + E_2$ is an ordered linear chain where $n \ge 0$, 
or (type $II_n$) $\Gamma_1 = 2(E+A_1+ \cdots + A_n) + A_{n+1} + A_{n+2}$
where $n \ge 1$ and both $E+A_1 + \cdots + A_n$ and $A_{n+1} + A_n + A_{n+2}$
are ordered linear chains, or (type $II_0$) $\Gamma_1 = A_1 + 2E + A_2$ 
where $A_1 + E + A_2$ is an ordered linear chain; here the $E, E_i$ are $(-1)$-curves 
and the $A_{\ell}$ are $(-2)$-curves.

\item[{\rm (4)}]
Let $\Delta=\sum_{i=1}^t\Delta_i$ be the decomposition into the connected components 
and let $(\Delta_i)_\red=\sum_{j=1}^{s_i}\Delta_i(j)$ be the irreducible decomposition 
with the dual graph below. We set $\Delta(j) := \Delta_1(j)$ when $t = 1$.

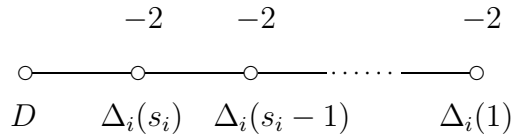
\raisebox{-25mm}{
\begin{picture}(75,25)(-20,0)
\unitlength=1mm
\put(5,15){\circle{1.8}}
\put(6,15){\line(1,0){13}}
\put(20,15){\circle{1.8}}
\put(21,15){\line(1,0){13}}
\put(35,15){\circle{1.8}}
\put(36,15){\line(1,0){9}}
\multiput(46,15)(1.5,0){6}{\circle*{0.2}}
\put(55,15){\line(1,0){9}}
\put(65,15){\circle{1.8}}
\put(3,8){$D$}
\put(15,8){$\Delta_i(s_i)$}
\put(30,8){$\Delta_i(s_i-1)$}
\put(60,8){$\Delta_i(1)$}
\put(18,21){$-2$}
\put(33,21){$-2$}
\put(63,21){$-2$}
\end{picture}}

\noindent
Then we have 
\[
\Delta_i=\sum_{j=1}^{s_i}\frac{j}{s_i+1}\Delta_i(j).
\]

\end{enumerate}
\end{lem}

\Proof
(1) Note that both $K_{\ol{X}}$ and $\ol{D}$ are in $(\Pic\ol{X})^G\otimes\Q$. 
We shall show that 
$-(K_{\ol{X}}+\ol{D})$ is $\Q$-ample. Note that $\kd(X,K_X+D+\Delta_{\red})=-\infty$. 
Suppose either $K_{\ol{X}}+\ol{D} \equiv 0$ or $K_{\ol{X}}+\ol{D}$ is $\Q$-ample. 
Consider the pull-back $K_X+D+\Delta=
f^*(K_{\ol{X}}+\ol{D})$. Then either $n(K_X+D+\Delta) \sim 0$ because $X$ is rational or 
$n(K_X+D+\Delta) >0$  for a positive integer $n$. Then we have 
\[
-\infty =\kd(X,K_X+D+\Delta_\red)\ge \kd(X,K_X+D+\Delta) \ge 0,
\]
which is absurd. So, $-(K_{\ol{X}}+\ol{D})$ is a $\Q$-ample divisor. Clearly, $\ol{D}$ is 
$\Q$-ample. Since $-K_{\ol{X}}=-(K_{\ol{X}}+\ol{D})+\ol{D}$, 
the divisor $-K_{\ol{X}}$ is ample.
\svskip

\noindent
(2) follows from the ampleness of $-K_{\ol X}$ and that $-K_X = f^*(-K_{\ol X})$.
\svskip

\noindent
(3) is a consequence of (2) (see also Lemma 1.3 in \cite{Z1}).
\svskip

\noindent
(4) follows from the fact that $-(K_X + D + \Delta) . \Delta_i(j) = 0$.
\QED

The following result describes roughly the shape of the divisor $D$.

\begin{lem}
We have the following assertions.
\begin{enumerate}
\item[{\rm (1)}]
Either $D \cong \BP^1$, or $D=D_1+D_2$, where $D_i \cong \BP^1$ and $D_1.D_2=1$. 
In both cases, $D.(D+K_X)=-2$.
\item[{\rm (2)}]
$D.\Delta < 2$, and $D_i.\Delta < 1$ if $D=D_1+D_2$.
\item[{\rm (3)}]
\[
0 < (K_X+D+\Delta)^2=(K_X+D).K_X+\Delta.D-2 < (K_X+D).K_X.
\]
\item[{\rm (4)}]
Suppose $X$ contains a $(-1)$-curve $E$ which is not a component of $D$. Then 
$E\cap D =\emptyset$ and $E.\Delta < 1$. Furthermore, $\Delta \ne 0$ and $E.\Delta > 0$. 
{\em See also Lemma 8.}
\end{enumerate}
\end{lem}

\Proof
Let $D_1$ be an irreducible component of $D$. Then we have 
\begin{eqnarray*}
\lefteqn{0 < -(K_X+D+\Delta).D_1 }\\
&&=2-2p_a(D_1)-(D-D_1).D_1-\Delta.D_1 \le 2-D_1.(D-D_1).
\end{eqnarray*}
This implies that $D_1.(D-D_1) \le 1$ and $D_1 \cong \BP^1$. Since $D+\Delta=f^*(\ol{D})$ 
is nef and big, it is $1$-connected, whence connected; see the proof of
Lemma 1 in \cite{Bo}. We note that if $D$ is not connected 
then $D+\Delta$ is not connected either. This is because each connected component of $\Delta$ 
meets exactly one irreducible component of $D$. The assertions (1) and (2) are thus proved. 
To verify the assertion (3), note that $\Delta.(K_X+D+\Delta)=0$ and $K_X.\Delta=0$. 
In view of the assertions (1) and (2), the computation is made as follows:
\begin{eqnarray*}
\lefteqn{0 < (K_X+D+\Delta)^2=(K_X+D+\Delta).(K_X+D)} \\
&& =(K_X+D)^2+\Delta.D=-2+(K_X+D).K_X+\Delta.D < (K_X+D).K_X.
\end{eqnarray*}
Let $E$ now be a $(-1)$-curve as in the assertion (4). Then we have
\[
0 < -(K_X+D+\Delta).E = 1-D.E-\Delta.E,
\]
where $\Delta.E \ge 0$. This implies that $D.E=0$ and $E.\Delta < 1$. 
Suppose $E.\Delta=0$. Then $E \cap (D + \Delta) = \emptyset$
and the image on ${\ol X}$ of $E$ is 
disjoint from $\ol{D}$. This contradicts the ampleness of $\ol{D}$. 
\QED

In case $\ol{X}$ has no singular points on $\ol{D}$, we can determine the surface $X$. 

\begin{lem} We have the following assertions.

\begin{enumerate}
\item[{\rm (1)}] If $K_X^2 \ge 8$ then one of the following cases occurs,
where $\Sigma_2$ is the Hirzebruch surface with a minimal
section $M$ such that $M^2 = -2$.
\item[{\rm (1a)}]
$X = {\BP^2}$ and $\deg D = 1, 2$, 
\item[{\rm (1b)}]
$X = {\BP^1} \times {\BP^1}$ and ${\rm bideg}\! D = (1, 1)$;
there is an element $g$ in $G$ which interchanges the
two different rulings on $X$,
\item[{\rm (1c)}]
${\ol X} = \ol{\Sigma}_2$ (the quadric cone in $\BP^3$)
and ${\ol D}$ is a hyperplane not passing through the vertex of the cone,
\item[{\rm (1d)}] 
$X= \Sigma_2$, $\Delta_{\rm red} = M$ and $D$ is a fibre,
\item[{\rm (1e)}] 
$X= \Sigma_2$, $\Delta_{\rm red} = M$ and $D$ is a section with self-intersection $4$,
\item[{\rm (1f)}] 
$X= \Sigma_2$, $\Delta_{\rm red} = M$ and $D$ is the sum of a fibre
and a section (disjoint from $M$) with self-intersection $2$.
\item[{\rm (2)}] If $\Delta = 0$ then $K_X^2 \ge 8$
and hence case (1a), (1b) or (1c) occurs.

\end{enumerate}
\end{lem}

\Proof
For the assertion (1), note that $X = \BP^2$ or $X$ is a Hirzebruch surface $\Sigma_d$
of degree $d \le 2$ (see Lemma 2). So if $\Delta \ne 0$, then $X = \Sigma_2$ and 
$\Delta = (1/2)M$. Now (1) follows from the fact that both $D + \Delta$ and 
$-(K_X + D + \Delta)$ are nef and big. The last part in $(1b)$ follows from the fact
that $\rho({\ol X}\quot G) = 1$.
\par

Let $\Delta = 0$. Suppose the contrary that
$K_X^2 \le 7$. Then $X$ contains a $(-1)$-curve $E$.
By Lemma 3, $E$ is contained in $D$. Since $D$ is nef and big
we have $D = D_1 + D_2$ with $D_1 = E$ and $D_2^2 \ge 0$.
Then both $D_i$ are $G$-stable. 
Hence $D_1=\alpha D_2$ with $\alpha > 0$. Indeed, $D_i$ is the total transform of its image on 
$\ol{X}\quot G$, and the images of the $D_i$ on $\ol{X}\quot G$ differ by a constant 
multiple which is a rational number $\alpha > 0$. Then we have 
\[
-1 = D_1^2=\alpha D_1.D_2=\alpha > 0,
\]
which is a contradiction.
Now (2) follows from (1).
\QED

{\em From now on, we assume that $\Delta \ne 0$.}

\begin{lem}
If $D=D_1+D_2$, we may assume that $D_2^2 \le 0$.
\end{lem}

\Proof
Suppose the contrary that $D_i^2 \ge 1$ for both $i=1,2$. We have 
\[
h^0(X,D_i)=D_i^2+2,
\]
which follows from an exact sequence
\[
0 \lto \SO_X \lto \SO_X(D_i) \lto \SO_{\BP^1}(D^2_i) \lto 0.
\]
Then one can find a member $\wt{D}_1 \in |D_1|$ such that $\#(\wt{D}_1\cap D_2)\ge D_1^2+1 
\ge 2$. Then $D_2$ is a component of $\wt{D}_1$ since $D_1.D_2=1$. Similarly, there exists 
a member $\wt{D}_2 \in |D_2|$ such that $D_1$ is a component of $\wt{D}_2$. This implies that 
$D_1 \sim D_2$ and $D_1^2=1$. Since $\Delta \ne 0$, there is
an irreducible component $\Delta_i(s)$ of $\Delta$ such that 
$D . \Delta_i(s) = 1$. This is absurd for $D = D_1 + D_2 \sim 2D_1$. 
So it is wrong to assume that $D_i^2\ge 1$ 
for both $i=1,2$. Hence we may assume that $D_2^2 \le 0$.
\QED

\begin{lem}
Suppose that both $D_1$ and $D_2$ are $G$-stable and $K_X^2 \le 7$. 
Then $\Delta$ has two connected components $\Delta_i$, where $\Delta_1$
is irreducible and $\Delta_2$ has length $2$. 
$f^{-1}{\ol D} + f^{-1}(\Sing {\ol X}) = D + \Delta_{\red} + A$
is given in Figure $1$ in Section 2.
\end{lem}

\Proof
We use the notation $\Delta = \sum_{i=1}^t \Delta_i$ and $\Delta_i = \sum_{j=1}^{s_i}
\Delta_i(j)$ of Lemma 2. We may assume that $\Delta_i$ meets $D_1$ (resp. $D_2$)
for $1 \le i \le t_1$ (resp. for $t_1 + 1 \le i \le t = t_1+t_2$).
By Lemma 3, $1 > D_1 . \Delta = \sum_{i=1}^{t_1} s_i/(s_i+1) \ge t_1/2$.
Hence $t_1 = 0, 1$. Similarly, $t_2 = 0, 1$. Thus $t = t_1 + t_2 = 1, 2$.
So we may assume that $D_i$ meets the connected component $\Delta_i$ of $\Delta$
of length $s_i$. We put $s_i = 0$ if $t_i = 0$.
\par

Since $\ol{D}_1=\alpha\ol{D}_2$ on $\ol{X}$ with 
$\alpha > 0$, we have 
$$
D_1+\Delta_1=\alpha(D_2+\Delta_2). \eqno{(1)}
$$
Taking the intersection of (1) with $D_1$, we have 
$$
D_1^2+\frac{s_1}{s_1+1}= \alpha. \eqno{(2)}
$$
Taking the intersection of (1) with $D_2$, we obtain
$$
1=\alpha\left(D_2^2+\frac{s_2}{s_2+1}\right). \eqno{(3)}
$$
Since $D_2^2$ is an integer, the equation (3) implies that $D_2^2 \ge 0$. Then $D_2^2=0$ 
by Lemma 5 and hence $\alpha=(s_2+1)/s_2$. Plugging the value of $\alpha$ in the equation (2), 
we have 
\[
D_1^2=\alpha-\frac{s_1}{s_1+1}=\frac{s_2+1}{s_2}-\frac{s_1}{s_1+1}=\frac{1}{s_2}+\frac{1}{s_1+1}.
\]
Since $D_1^2$ is an integer, we have the following possibilities:
\[
(s_1,s_2;\alpha;D^2_1) = (0,1;2;2), \, (1,2;3/2;1).
\]
We shall show the assertion that the first case (resp. the second case) implies that 
$X = \Sigma_2$ and $K_X^2 = 8$ (resp. implies the result of the lemma).
Indeed, in the first (resp. second) case, if $\Gamma_2$ is a singular fibre of the 
$\BP^1$-fibration $\varphi$ induced by $|D_2|$ containing no components of $\Delta$ 
(resp. containing $\Delta_2(1)$ but no $\Delta_1$), then it is of type $I_n$ in Lemma 2 and 
the cross-section $D_1$ must meet a $(-1)$-curve $E_1$ 
of $\Gamma_2$, contradicting Lemma 3 (4). 
In the second case, if $\Gamma_2$ is the fibre of $\varphi$ containing both
$\Delta_2(1)$ and $\Delta_1$ then it is of type $II_0$.
So the assertion is true and the lemma proved.
\QED

Now consider the case where $D=D_1+D_2$ and $g(D_1)=D_2$ for some $g$ in $G$. 
Then $D_1^2 =D_2^2$.  By the proof of Lemma 6, $\Delta$ has two connected
components $\Delta_i$ and we may assume that
$D_i$ meets $\Delta_i$ so that $g(\Delta_1) = \Delta_2$ and hence
$\Delta_1$ and $\Delta_2$ are $(-2)$-linear chains of the same length $s \ge 1$.
Since $\ol{D}$ is not contractible to a point, it follows that $D^2_1=D_2^2 \ge -1$. By 
Lemma 5, we have $D_1^2=-1$ or $0$. Write 
$(\Delta_i)_{\rm red} = \sum_{j=1}^{s_i} \Delta_i(j)$ as in Lemma 2.

\begin{lem}
Suppose that $g(D_1)=D_2$ for some $g$ in $G$. Then $D_1^2=-1$. 
Furthermore, there is a $\BP^1$-fibration $\Phi : X \rightarrow B$ ($B \cong {\BP^1}$)
such that $D_1+D_2$ is a fibre. According to the 
possible types of the singular fibres of $\Phi$, we have five different cases; 
see Figures $2 \sim 6$ in Section 2,
each of which also contains the graph 
$f^{-1}{\ol D} + f^{-1}(\Sing {\ol X}) = D + \Delta_{\red} + A$.
\par
$G = \langle g \rangle \cong {\Z}/(2)$ is realizable; indeed
either $X$ or its blow-down of the $G$-stable curve $E$ (for Figures $2$ and $3$)
is obtained from the Hirzebruch surface 
$\Sigma_2$ (or one point $Q_i$ blow-up of $\Sigma_2$ for Figures $2, 5, 6$)
by taking a double cover ramifying along a smooth
(or singular at $Q_i$) irreducible member of 
$|-K_{\Sigma_2}|$ with $G$ equal to
the Galois group ${\rm Gal}(X/\Sigma_2)$. 
\end{lem}

\Proof
Suppose $D_1^2 \ne -1$. Then $D_1^2=D_2^2=0$ as shown above. Consider the $\BP^1$-fibraton 
$\Phi_{|D_2|}$ for which $D_1$ and the component $\Delta_2(s)$ are cross-sections. 
Since $\Delta_1 \ne 0$, the map $\Phi_{|D_2|}$ has a singular fibre $\Gamma_1$ comprising 
$(\Delta_1)_\red$ and $(\Delta_2)_\red-\Delta_2(s)$. Note that there are no other 
singular fibres because no $(-1)$-curves lying outside of $D$ meet the cross-section 
$D_1$ by Lemma 3. The singular fibre $\Gamma_1$ consists of $(-2)$-curves and $(-1)$-curves. 
Since $\Delta_1$ and $\Delta_2$ have the same length $s$, only possibility for 
$\Gamma_1$ is that $s=1$ and the dual graph of $\Gamma = E_1 + \Delta_1(1) + E_2$ is 
of type $I_1$ in Lemma 2 such that $E_1$ meets the cross-section $\Delta_2(1)$.
Then $E_1.\Delta=1$, which is a contradiction to Lemma 3. 
\svskip 

Now assume that $D_1^2=D_2^2=-1$. Then $|D_1+D_2|$ defines a $\BP^1$-fibration $\Phi :
X \to B$, where $B \cong \BP^1$. 
The $\Delta_i(s)$ are the cross-sections of $\Phi$. 
Suppose that $\Delta_1$ has length $s \ge 2$. Let $\Gamma_1$ be the singular fibre of $\Phi$ 
containing $(\Delta_1)_\red-\Delta_1(s)=\Delta_1(s-1)+\Delta_1(s-2)+\cdots+\Delta_1(1)$. 
If $\Gamma_1$ contains no components of $\Delta_2$, then $\Gamma_1 = 
E_1 + \Delta_1(s-1) + \cdots +\Delta_1(1) + E_2$
is of type $I_{s-1}$ in Lemma 2 so that $E_i$ meets the cross-section $\Delta_2(s)$
for $i = 1$ or $2$. Then we have by Lemma 2,
\[
\Delta.E_i \ge \frac{1}{s+1}+\frac{s}{s+1} = 1,
\]
which is impossible by Lemma 3. 
So, the fibre $\Gamma_1$ contains also $(\Delta_2)_\red-\Delta_2(s)$.
This implies that $s=2$ and $\Gamma_1=\Delta_1(1)+2E+\Delta_2(1)$ is of
type $II_0$ in Lemma 2, where $E$ is a $(-1)$-curve. In 
particular, the length of the $(-2)$ linear chain $\Delta_i$ is less than or equal to $2$. 
The possible cases of all singular fibres of $\Phi$ are exhausted 
by the following four (see Fact 12 in the proof of Lemma 11);
this proves the lemma (the realization part is easy to check).

\begin{enumerate}
\item
$s=2$; $K_X^2 = 3$; $\Gamma_0:=D_1+D_2, \Gamma_1:=\Delta_1(1)+2E+\Delta_2(1)$ and 
$\Gamma_2:= E_1+A+E_2$ are of types $I_0, II_0, I_1$ in Lemma 2. See Figure 2.
\item
$s=2$; $K_X^2 = 3$; $\Gamma_0:=D_1+D_2, \Gamma_1:=\Delta_1(1)+2E+\Delta_2(1)$ and
$\Gamma_{i+1}:= E_i + F_i$ ($i = 1, 2$) are of types $I_0, II_0, I_0, I_0$.
See Figure 3.
\item
$s=1$; $K_X^2 = 4$; $\Gamma_0:= D_1+D_2$ and $\Gamma_i:=E_i+F_i$ ($i=1,2,3$) are all of types
$I_0$. See Figure 4.
\item
$s=1$; $K_X^2 = 4$; $\Gamma_0:=D_1+D_2$ and $\Gamma_1=E_1+A_1+A_2+E_2$ are of types $I_0, I_2$.
See Figure 5.
\item
$s=1$; $K_X^2 = 4$; $\Gamma_0:=D_1+D_2$, $\Gamma_1=E_1+F_1$
and $\Gamma_2 = E_2 + A + F_2$ are of types $I_0, I_0, I_1$.
See Figure 6.
\end{enumerate}
\QED

Now we switch to the case $D$ is irreducible. Note that $\Delta\ne 0$ is assumed.

\begin{lem}
Suppose $D$ is irreducible. Then the following assertions hold.
\begin{enumerate}
\item[{\rm (1)}] 
$h^0(-K_X-D) \ge K_X.(K_X+D)>0$. Hence $|-(K_X+D)|\ne \emptyset$.
\item[{\rm (2)}]
$|-(K_X+D+\Delta_\red)|\ne \emptyset$.
\item[{\rm (3)}]
Let $E$ ($\ne D$) be a $(-1)$-curve. Then $E.\Delta_\red=1,2$. 
If $E.\Delta_\red=2$, then $D.\Delta_\red=2$ (i.e., $t = 2$ 
in notation of Lemma $2$) and
$E+D+\Delta_\red$ is a simple loop and linearly equivalent to $-K_X$.
\item[{\rm (4)}]
We have $D^2 \ge -1$. 
The number $t$ of connected components of $\Delta$ is at most $3$.
If $t = 3$, then, in notation of Lemma $2$,
$(s_1+1, s_2+1, s_3+1) = (2, 2, n)$ $(n \ge 2)$ or $(2, 3, n)$
$(n = 3, 4, 5)$ {\em (those triplets are called the Platonic numbers)}.
\end{enumerate}
\end{lem}

\Proof
By the Riemann-Roch theorem, we have 
\begin{eqnarray*}
\lefteqn{h^0(-K_X-D)-h^1(-K_X-D)+h^0(2K_X+D)} \\
&& = \frac{(-K_X-D).(-2K_X-D)}{2}+1 \\
&& =K_X^2+\frac{D^2}{2}+\frac{3K_X.D}{2}+1\\
&& =K_X^2+K_X.D >0;
\end{eqnarray*}
here we used that $K_X.D+D^2=-2$ and $K_X^2+K_X.D > 0$ in Lemma 3. 
Note that $h^0(2K_X+D) \le h^0(2(K_X+D+\Delta))=0$
because $-(K_X+D+\Delta)$ is nef and big. Now the assertion (1) follows.
\par

Let $({\Delta_1})_\red=A_s+A_{s-1}+\cdots+A_1$ be a connected component of $\Delta_\red$ such that 
$D.A_s=1$. Since $-(K_X+D).A_s=-D.A_s <0$, we have $|-(K_X+D+A_s)|\ne \emptyset$. Suppose 
$|-(K_X+D+A_s+\cdots+A_i)|\ne \emptyset$. Since $-(K_X+D+A_s+\cdots+A_i).A_{i-1}=-A_i.A_{i-1}
< 0$, it follows that $|-(K_X+D+A_s+\cdots+A_{i-1})|\ne \emptyset$. So
$|-(K_X+D+({\Delta_1})_\red)|\ne \emptyset$. Likewise, $|-(K_X+D+\Delta_\red)|\ne \emptyset$.
\par

Let $E$ be a $(-1)$-curve not in $D$. Then $E.\Delta_\red > 0$ by Lemma 3. Let $p : X \to Y$ be the 
blow-down of $E$. If $E.\Delta_\red\ge 2$ then $|K_Y+p_*(D+\Delta_\red)|\ne \emptyset$ by 
the Riemann-Roch theorem or Lemma 2.1.3 in \cite{M1}, page 7. 
Since $|-(K_Y+p_*(D+\Delta_\red))|\ne \emptyset$ by the assertion (2), it 
follows that $K_Y + p_*(D + \Delta_{\red}) \sim 0$.
So $t = 2$ and $E$ meets the end
component of each $\Delta_i$ which is located on the opposite side of $D$. 
Namely, $D+\Delta_\red+E$ is a simple loop. 

\par
Since $D.\Delta < 2$ by Lemma 3, in notation of Lemma 2, we have 
\[
\sum_{i=1}^t\left(1-\frac{1}{s_i+1}\right) < 2.
\]
It then follows that $t \le 3$. Moreover, if $t = 3$ then $\{s_1+1,s_2+1,s_3+1\}$ 
is one of the Platonic triplets upto permutations. So always $D^2 \ge -1$, for otherwise
$D + \Delta$ is contractible, contradicting the nef and bigness of $D + \Delta$
(see Satz 2.11 in \cite{Br}).
\QED

Consider the case $\Delta$ has three connected components $\Delta_i$
(see Lemma 2). Let $C_i$ be the component of $\Delta_i$ meeting $D$,
i.e., $C_i = \Delta_i(s_i)$ in notation of Lemma 2.

\begin{lem}
Suppose that $D$ is irreducible and $\Delta$ has three connected components.
Then $D^2=-1$. Also $\Delta_\red$ consists of three disjoint irreducible 
curves $C_i=(\Delta_i)_\red \ (i=1,2,3)$ and $-K_X= 2D+C_1+C_2+C_3$. Furthermore, $K_X^2=2$ and 
there is a birational morphism $q : X \to \BP^2$ such that $q_*(D+C_1+C_2+C_3)$ is a union 
of a line and a conic touching each other in one point.
$|2D+C_1+C_2|$ defines a $\BP^1$-fibration $\Phi : X \to B$, 
and according to the different 
types of the singular fibres of $\Phi$, there are seven possible cases;
see Figures $7 \sim 13$, each of which also contains the graph 
$f^{-1}{\ol D} + f^{-1}(\Sing {\ol X}) = D + \Delta_{\red} + A$.
\end{lem}

\Proof
By Lemma 8, $D^2 \ge -1$. Consider the case $D^2 = -1$.
Let $p : X \to Y$ be the blow-down of $D$ and let $\ol{C}_i=p(C_i)$. Then the $\ol{C}_i$ 
share one point in common. So $|K_Y+\ol{C}_1+\ol{C}_2+\ol{C}_3|\ne \emptyset$ 
by the Riemann-Roch theorem or Lemma 2.1.3 in \cite{M1}, page 7.
Since $|-(K_Y+\ol{\Delta}_\red)|\ne \emptyset$ by Lemma 8, we have $-(\ol{\Delta}_\red-
\ol{C}_1-\ol{C}_2-\ol{C}_2) \ge 0$, where $\ol{\Delta}=p_*(\Delta)$. Hence it follows that 
$\ol{\Delta}_\red=\ol{C}_1+\ol{C}_2+\ol{C}_3$ and $K_Y+\ol{C}_1+\ol{C}_2+\ol{C}_3=0$. Thence 
follows the assertion on the expression of $-K_X$. In order to obtain the morphism $q$,
we let $q_1 : X \to X_1$ be the blow-down of $D + C_3$ and 
continue blowing down further to reach a relatively minimal model $\Sigma_d$
with $d = 0,1,2$ (see Lemma 2). Then one can bypass the blow-down steps to reach $\BP^2$. 
By making use of the property that $-K_Z$ for a surface $Z$ appearing in the blow-down step 
is the sum of the images of $C_1$ and $C_2$, any $(-1)$-curve on $Z$ meets transversally 
exactly one of the images of $C_1$ and $C_2$ in one point. If we set $B_1:=q(C_1)$ and 
$B_2:=q(C_2)$, then $B_1+B_2$ is a cubic curve and $B_1\cap B_2$ consists of a single point. 
Hence we may assume that $B_1$ is a line 
and $B_2$ is a conic. Let $\Gamma_0 := 2D+(\Delta_1)_{\red}+(\Delta_2)_{\red}$ 
and let $\Phi$ be the ${\BP^1}$-fibration with $\Gamma_0$ as a fibre.
Since $-K_X = 2D+\sum_i (\Delta_i)_{\red}$ supports a fibre and
a 2-section $(\Delta_3)_{\red}$, every $(-2)$-curve, i.e., every
component of $f^{-1}(\Sing {\ol X})$ other than $(\Delta_3)_{\red}$ 
is contained in a fibre.
So $f^{-1}({\ol D}) + f^{-1}(\Sing {\ol X}) = D + \Delta_{\red} + A$
is given in one of Figures $7 \sim 13$ in Section 2. 
See Lemma 2 for possible types of singular fibres;
see also Fact 12 in the proof of Lemma 11.
To be precise, the following cases are not included but reduced to other cases,
and Figures 7-7' appear on the same $X$ with two different fibrations.

\par \svskip \noindent
{\sc Case 9.1.}\ \ $\Gamma_0, \Gamma_1 = E_1 + A_1 + A_2 + A_3 + E_2$ which is
of type $I_3$ in Lemma 2, are the only singular fibres of $\Phi$. 
Also the 2-section $(\Delta_3)_{\red}$ meets each $E_i$. By going to
a Hirzebruch surface $\Sigma_d$ ($d \le 2$), we see that there is a $(-1)$-curve
$E$ on $X$ such that $E . A_2 = E . (\Delta_i)_{\red} = 1$ for $i = 1$ or $2$
say for $i = 1$. Then $\Gamma_0' := 2D + (\Delta_2)_{\red} + (\Delta_3)_{\red}$ is the
singular fibre of a new ${\BP^1}$-fibration $\Phi'$, and
$\Gamma_1' := 2(E + A_2) + A_1 + A_3$ is also a singular fibre of $\Phi'$.
So $\Phi'$ fits Figure 8 after relabeling $\Delta_i$.
\svskip 

\noindent
{\sc Case 9.2.}\ \ $\Gamma_0$, and $\Gamma_i = E_i + A_i + F_i$ ($i = 1, 2$)
each of which is of type $I_1$ in Lemma 2, are the only singular fibres of $\Phi$. 
Also the 2-section $(\Delta_3)_{\red}$ meets each of $E_i, F_j$. 
We can find a $(-1)$-curve $E$ on $X$ such that 
$E . A_1 = E . A_2 = E . \Delta_i = 1$ for $i = 1$ or $2$
say for $i = 1$. Then $\Gamma_0' := 2D + (\Delta_2)_{\red} + (\Delta_3)_{\red}$ is the
singular fibre of a new ${\BP^1}$-fibration $\Phi'$, and
$\Gamma_1' := 2E + A_1 + A_2$ is also a singular fibre of $\Phi'$.
So $\Phi'$ fits Figure 10 after relabeling $\Delta_i$.
\svskip 

Consider the case $D^2=0$. Then $|D|$ defines a $\BP^1$-fibration $\Phi : X \to B$ 
for which the curves $C_1, C_2, C_3$ are cross-sections. 
Suppose $\{s_1+1,s_2+1,s_3+1\}=\{2,2,n\}$ with $n \ge 3$. Write $\Delta_3=C_3+A_m+\cdots+A_1$ 
with $m=n-2 \ge 1$. Then there exists a singular fibre $\Gamma_1$ of $\Phi$ such that 
$\Gamma_1=E_1+A_m+\cdots+A_1+E_2$ is an ordered linear chain and of type $I_m$ in Lemma 2
so that $E_i$ meets the cross-section $C_2$ for $i = 1$ or $2$.
Then $E_i + D + \Delta_{\red}$ contains a loop and $(\Delta_1)_{\red}$,
contradicting Lemma 8. In the case $\{s_1+1,s_2+1,s_3+1\}=\{2,2,2\}$, we note that 
$\Phi$ is not a relatively minimal $\BP^1$-fibration. Hence $\Phi$ has a singular fibre 
$\Gamma_1$ of type $I_k$ and $\Gamma_1$ contains a $(-1)$-curve $E_1$ meeting two of 
$C_1,C_2,C_3$, say meeting $C_1, C_2$. Then $E+D+\Delta_\red$ contains a loop and $C_3$, 
contradicting Lemma 8. The case $\{s_1+1,s_2+1,s_3+1\}=\{2,3,n\}$ with $n = 3,4,5$ also 
leads to a contradiction.
\svskip 

Consider the case $D^2=1$. Let $p : X \to \BP^2$ be the birational morphism defined by $|D|$.
Since $D^2=1$ and hence $D$ is linearly equivalent to the pull-back of a line, the morphism 
$p$ is a composite of the blow-downs of $(-1)$-curves which are disjoint from $D$ and 
its images. This implies that $B_i:=p(C_i)\ (i=1,2,3)$ is a curve. Since 
\[
-K_{\BP^2}=-p_*(K_X) \ge p_*(D+C_1+C_2+C_3),
\]
it follows that $\deg(-K_{\BP^2}) \ge 4$, which is a contradiction.
\svskip 

Consider the case $D^2 \ge 2$. Then we have by Lemma 3
\[
K_X^2 \ge 1-D.K_X=3+D^2 \ge 5.
\]
Hence we have 
\[
5 \ge 10-K_X^2=\rho(X) \ge \rho(\ol{X})+\#\Delta \ge 1+\#\Delta,
\]
where $\#\Delta$ signifies the number of the irreducible components of $\Delta$. 
So, $s_1 + s_2 + s_3 = \#\Delta \le 4$. Thus the possible 
cases of $\{s_1+1,s_2+1,s_3+1\}$ are $\{2,2,2\}$ and $\{2,2,3\}$ up to permutations. 
If $\{s_1+1,s_2+1,s_3+1\}=\{2,2,3\}$, then $\rho(\ol{X})=1, \Sing\ol{X}=2A_1+A_2$ and 
$K_X^2=5$. But this case cannot occur by the classifications of the distributions of 
singular points (cf. Lemma 3 in part I of \cite{MZ}). If $\{s_1+1,s_2+1,s_3+1\}=\{2,2,2\}$, 
then either $\rho(X)=4$ or $\rho(X)=5$. In the first case, we have $\rho(\ol{X})=1$ and 
$\Sing\ol{X}=3A_1$, which is also impossible \cite{MZ}. 
In the second case, we have either $\rho(\ol{X})=1$ and $\Sing\ol{X}=4A_1$, or 
$\rho(\ol{X})=2$ and $\Sing\ol{X}=3A_1$. The case $\rho(\ol{X})=1$ is ruled out by 
\cite{MZ} and the case $\rho(\ol{X})=2$ by \cite{Ye}.
\QED

Next we consider the case where $\Delta$ is connected. 
Write $\Delta_{\red} = \Delta(s) + \cdots + \Delta(1)$ as an ordered linear chain
so that $D . \Delta(s) = 1$.

\begin{lem}
Suppose that $D$ is irreducible,
$\Delta$ is connected of length $s$ and $K_X^2 \le 7$. 
Then the following assertions hold, where each of Figures $14 \sim 23$
contains the graph $f^{-1}{\ol D} + f^{-1}(\Sing {\ol X}) = D + \Delta_{\red} + A$.
\begin{enumerate}
\item[{\rm (1)}]
For any $(-1)$-curve $E$ on $X$, it holds that $E.\Delta_\red=1$ and $E\cap D=\emptyset$.
\item[{\rm (2)}]
$D^2=0,1,2$.
\item[{\rm (3)}]
If $D^2=0$, then $s = 2$ (resp. $4$) and $K_X^2=6$ (resp. $K_X^2=5$). There are three 
possible cases (see Figures $14 \sim 16$).
In Figure $16$, there is an element $g$ in $G$ such that $g(E_1) = E_2$.
\item[{\rm (4)}]
If $D^2=2$, then $s=2, K_X^2=6$ and there are two possible cases 
(see Figures $17$ and $18$). In Figure $17$ (resp. $18$), $E$ (resp. $E_1, E_2$)
are the only $(-1)$-curve(s) on $X$.
In Figure $18$, we have $g(E_1) = E_2$ for some $g$ in $G$.

\item[{\rm (5)}]
If $D^2=1$, then either $s=4$ and $K_X^2 = 5$ (see Figure $19$), 
or $s=1$ and $K_X^2 = 6$ (see Figures $20-21$).
In Figures $19-20$, the $E$ is the only $(-1)$-curve on $X$.
\end{enumerate}
\end{lem}

\Proof
(1)\ It follows from Lemmas 3 and 8 (see also (2)).
\svskip

\noindent
(2)\ By Lemma 3, we have
\[
7 \ge K_X^2 > 2-\Delta.D-K_X.D=4+D^2-\Delta.D=3+D^2+\frac{1}{s+1}.
\]
Since $D+\Delta$ is nef and big, we have $D^2 \ge 0$, for otherwise $D+\Delta_\red$ is 
contractible to a point, a contradiction. 
Hence $7 \ge K_X^2 \ge 4+D^2 \ge 4$. So, $D^2 = 0,1,2,3$. 
\par

Suppose $D^2=3$. Then $K_X^2=7$ and 
\[
3=\rho(X)\ge \rho(\ol{X})+s\ge 1+s \ge 2,
\]
whence $s=1,2$. If $s=2$, then $\rho(\ol{X})=1$ and $\Sing\ol{X}=A_2$, and this case 
does not occur (cf. \cite{MZ}). So, $s=1$. If $\rho(\ol{X})=1$, then $\Sing\ol{X}=A_1+A_1$, 
and this case does not occur either (cf. {\em ibid.}). Thus $\rho(\ol{X})=2, K_X^2=7$ and 
$\Sing\ol{X}=A_1$. 
\par

Take a $(-1)$-curve $E_1$ on $X$. Then $E_1.\Delta_\red=1$ by the assertion (1). 
The blow-down of $E_1$ brings $X$ to the Hirzebruch surface $\Sigma_1$, 
as the image of $\Delta_\red$ becomes a $(-1)$-curve. Let $P$ be the image of $E_1$ on $\Sigma_1$. 
Then the proper transform $E_2$ of a fibre of the ruling on $\Sigma_1$ 
passing through $P$ is a $(-1)$-curve meeting $D$. 
This contradicts the assertion (1). Hence $D^2=0,1,2$.
\svskip 

\noindent
(3)\ Let $D^2 = 0$. Let $\Phi : X \to B$ be the $\BP^1$-fibration with $D$
as a fibre. Since $D$ is $G$-stable, $G$ permutes fibres of $\Phi$.
\par

Suppose $D^2=0$ and $s=1$. Then there is a singular fibre 
$\Gamma=E_1+A_1+\cdots+A_n+E_2$ of type $I_n$, where we may assume that
the $(-1)$-curve $E_1$ meets the cross-section $\Delta_\red$. 
Then $E_2 \cap \Delta = \emptyset$, contradicting Lemma 3.
\par

Suppose $D^2=0$ and $s \ge 2$. Let $\Gamma_1$ be the singular fibre containing 
$\Delta(s-1) + \cdots + \Delta(1)$. Then $\Gamma_1$ is $G$-stable.
If $\Phi$ contains a second singular fibre $\Gamma_2$, then we can reach the same 
contradiction as above. Hence $\Gamma_1$ is the only singular fibre of $\Phi$.
If $s \ge 3$ and $\Gamma_1 = E_1 + \Delta(1) + \cdots + \Delta(s-1) + E_2$
is an ordered linear chain and a singular fiber of type $I_{s-1}$, then the image on 
${\ol X}$ of $E_1$ is $G$-stable and contractible, contradicting $\rho({\ol X}\quot G) = 1$.
By the arguments above, all possible types of $\Gamma_1$ are given in
Figures $14 \sim 16$.  In Figure 16, since $\Gamma_1$ is $G$-stable, each
element in $G$ either stabilizes $E_i$ or interchanges $E_1, E_2$.
If $E_i$ is $G$-stable then the image on ${\ol X}$ of $E_i$ is $G$-stable
and contractible, contradicting $\rho({\ol X}\quot G) = 1$.
\svskip 

\noindent
(4)\ Suppose $D^2=2$. Since $7\ge K_X^2\ge 4+D^2=6$, we have $K_X^2=6$ or $7$. Accordingly, 
$\rho(X)=4$ or $3$. Since $2 \le 1+s \le \rho(\ol{X})+s \le \rho(X)$, we have $s=3, 2, 1$. 
\par

Suppose $s=3$. Then $\rho(X)=4, \rho(\ol{X})=1$ and $\Sing\ol{X}=A_3$. This case does not 
occur by \cite{MZ}. 
\par

Suppose $s=2$. If $\rho(X)=3$ then $\rho(\ol{X})=1$ and $\Sing\ol{X}=A_2$, 
and this case does not occur either \cite{MZ}. If $\rho(X)=4$, 
either $\rho(\ol{X})=1$, $\Sing\ol{X}
=A_1+A_2$ and there is only one $(-1)$-curve on $X$, 
or $\rho(\ol{X})=2$, $\Sing\ol{X}=A_2$ and there are only two
$(-1)$-curves on $X$; see Figures 5 and 6 in \cite{Ye}.
\par

Take a $(-1)$-curve $E_1$ on $X$. Note that $E_1.\Delta_\red=1$ by the assertion (1). 
Let $X \rightarrow {\BP^2}$ be the blow-down of $E_1 + \Delta(1) + \Delta(2)$
and ${\widetilde D}$ the image of $D$.
If $E_1 . \Delta(1) = 1$, then ${\widetilde D}^2 = 3$, which
is impossible on ${\BP^2}$. Hence $E_1 . \Delta(2) = 1$ and ${\widetilde D}^2 = 4$.
Let $P$ be the point on $\wt{D}$ which is the fundamental point of the blow-down $X \to \BP^2$. 
Let $\ell$ be a line which is tangent to $\wt{D}$ at $P$. Reverse the above blow-down. Let 
$L$ be the proper transform of $\ell$. There are two possibilities according as 
$E_1\cap L \ne \emptyset$ or $E_1\cap L= \emptyset$. In the former case, $L$ is a $(-2)$-curve 
(see Figure 17 where $A : = L$ and $\rho({\ol X}) = 1$) and, in the latter case, 
$L$ is a $(-1)$-curve (see Figure 18 where $E_2 := L$ and $\rho({\ol X}) = 2$).
Note that $E$ (resp. $E_1$ and $E_2$) is/are the only $(-1)$-curve(s) on $X$
(see Figures 5 and 6 in \cite{Ye}). 
In Figure 18, we have $g(E_1) = E_2$ for some $g$ in $G$ (see the argument for
Figure 16).
\par

Suppose $s=1$. Take a $(-1)$-curve $E_1$. Then $E_1 . \Delta(1) = 1$. 
Let $p : X \to Y$ be the blow-down of $E_1$ and 
$\Delta(1)$. Since $K_X^2 \ge 6$, we have $K_Y^2 \ge 8$ and $\wt{D}^2=3$. 
Since there are no curves $\wt{D}$ on $\BP^2$ with $\wt{D}^2=3$, $K_Y^2=8$. 
Hence $Y \cong \Sigma_d$ with $d=0,1,2$ (see Lemma 2).
But any curve on $\Sigma_d$ with $d=0, 2$ has self-intersection number divisible by $2$. So,
$Y\cong \Sigma_1$. Let $M$ and $\ell$ be respectively the minimal section and a fibre on 
$\Sigma_1$. Then $\wt{D}\sim M+2\ell$. Let $P$ be the fundamental point of $p$. If $P=
\wt{D}\cap M$, then $\Delta(1)+p'(M)$ is a $(-2)$-chain, for otherwise 
there appears a $(-3)$-curve on $X$; since $s=1$, this is a contradiction. 
If $P \ne \wt{D}\cap M$, then $p'(M)$ is a 
$(-1)$-curve meeting $D$, contradicting Lemma 3. 

\svskip \noindent
(5)\ Now assume that $D^2=1$. 
Note that $7 \ge K_X^2 \ge 4+D^2=5$. On the other hand, since 
\[
2 \le s+1 \le s+\rho(\ol{X}) \le \rho(X) = 10-K_X^2 \le 6-D^2=5,
\]
we have $1 \le s \le 4$. We consider all possible cases according to the value of $s$. 
We note that $E.\Delta_\red=1$ for any $(-1)$-curve $E$ on $X$ by the assertion (1).

\svskip \noindent
{\sc Case $s=4$.}\ \ Then $K_X^2=5$, $\rho({\ol X}) = 1$ and $\Sing {\ol X} = A_4$.
Take a $(-1)$-curve $E$ on $X$. If $E.\Delta(i)=1$ for 
$i=1$ or $4$, then one can blow down $E+\Delta_\red$ and the resulting surface $Y$ has 
$K_Y^2=10$. This is a contradiction. Suppose $E.\Delta(2)=1$. Then $D$ is a component of a 
fibre of the $\BP^1$-fibration defined by $|2(E+\Delta(2))+\Delta(1)+\Delta(3)|$. Since 
$D^2=1$, this is a contradiction. Consequently, $E.\Delta(3)=1$. 
We thus obtain Figure 19, where $E$ is the only $(-1)$-curve on $X$
(see Figure 5 in \cite{Ye}).

\svskip \noindent
{\sc Case $s=3$.}\ \ We claim that this case does not take place.
Note that $\rho(X)=4$ or $5$. If $\rho(X)=4$ then $\rho(\ol{X})=1$ and 
$\Sing\ol{X}=A_3$, which is not the case by \cite{MZ}. So $\rho(X)=5$. If $\rho(\ol{X})=1$ 
then $\Sing\ol{X}=A_1+A_3$, which is not the case either by \cite{MZ}. So $\rho(\ol{X})=2, \,
\Sing\ol{X}=A_3$ and $K_X^2=5$. Then there are only two $(-1)$-curves $E_1, E_2$ on $X$
with $E_i . \Delta_{\red} = 1$ ($i = 1, 2$) and $E_1 . \Delta(2) = E_2 . (\Delta(1) + \Delta(3)) = 1$
(see Figure 6 in \cite{Ye}).
So both $E_i$ are $G$-stable, and the image on ${\ol X}$ of $E_2$
is $G$-stable and contractible, contradicting $\rho({\ol X} \quot G) = 1$.

\svskip \noindent
{\sc Case $s=2$.}\ \ We shall show that this case does not take place. 
Let $E$ be a $(-1)$-curve. Then $E . \Delta_{\red} = 1$.
Let $p : X \to Y$ be the blow-down of $E+\Delta(1)+\Delta(2)$. Since $K_X^2 \ge 5$,
we have $K_Y^2 \ge 8$.
If $E . \Delta(2) = 1$, then $p_*(D)^2=3$. As in the proof of the assertion (4) for 
the case $D^2=2$ and $s=1$, $Y$ must be the Hirzebruch surface $\Sigma_1$. Let $M$ be the 
minimal section of $\Sigma_1$. Since $s=2$, the fundamental point of $p$ is different from 
the point $M\cap p_*(D)$. Then $E_1:=p'(M)$ is a $(-1)$-curve such that $E_1.D=1$, 
contradicting Lemma 3. Hence $E.\Delta(1)=1$. 
Then $p_*(D)^2=2$ and $K_Y^2\ge 8$. Hence $Y$ is the Hirzebruch surface 
$\Sigma_d$. Since $p_*(D)^2=2$, one can readily show that $d=0, 2$ (see Lemma 2). 
\par

Suppose first that $Y \cong \Sigma_2$. Let $M$ be the minimal section. 
Then $M\cap p_*(D)=\emptyset$. 
Reversing the above blow-down and noting that the length of $\Delta_\red$ is $2$, we 
can show that $E_1:=p'(\ell)$ is a $(-1)$-curve meeting $\Delta(2)$, where $\ell$ is 
the fibre passing through the fundamental point of $p$.  So, we are lead to a contradiction 
by the above case. 
\par

Suppose $Y \cong \Sigma_0$. Let $\ell$ be one of the fibres (of the 
two different $\BP^1$-fibrations) passing through the fundamental point of $p$. Then it 
follows that $E_1:=p'(\ell)$ is a $(-1)$-curve meeting $\Delta(2)$. Again, we are lead 
to a contradiction. Consequently, the case $s=2$ does not occur. 

\svskip \noindent
{\sc Case $s=1$.}\ \ Let $E$ be a $(-1)$-curve. Then $E.\Delta(1)=1$. Let $p : X \to Y$ 
be the blow-down of $E+\Delta(1)$. Since $K_X^2 \ge 5$, we have $K_Y^2 \ge 7$. 
\par

Suppose $K_X^2=5$. Then $Y$ has a $\BP^1$-fibration $\pi$, which is not relatively minimal but 
contains a singular fibre consisting of two $(-1)$-curves $E_1+E_2$. Since $p_*(D)^2=2$, 
$p_*(D)$ is not contained in a fibre of $\pi$. This implies that 
$p_*(D)\cap E_i \ne \emptyset$ for $i = 1$ or $i = 2$, say for $i = 1$.
Then the fundamental point $P$ of the morphism $p$ is not contained in $p_*(D)\cap E_1$, 
for otherwise $s \ge 2$ or $p'(E_1)^2 \le -3$, a contradiction. Hence $p'(E_1)$ remains 
as a $(-1)$-curve on $X$ which meets $D$.  This contradicts Lemma 3.
We have therefore $K_X^2=6$ and $Y \cong \Sigma_d$ with $d=0,2$ because $p_*(D)^2=2$
(see Lemma 2). 
\par

Suppose $Y \cong \Sigma_2$. Let 
$M$ be the minimal section and let $\ell$ be the fibre passing through the fundamental 
point $P$ of $p$. Consider the inverse of the morphism $p$. After blowing up the 
point $P$, there are two possibilities of taking the centre $Q$ of the next blow-up. 
Namely, $Q$ lies (resp. does not lie) on the proper transform of $\ell$. 
The first case gives rise to Figure 20, where $\rho({\ol X}) = 1$,
$\Sing {\ol X} = A_1 + A_2$ and $E$ is the only $(-1)$-curve on $X$
(see \cite{Ye}, Figure 5).
In the second case, $E + \Delta(1) + \ell' + M'$ has the dual graph
$$(-1)-(-2)-(-1)-(-2)$$
where $\ell', M'$ are the proper transforms of $\ell, M$.
If $\rho({\ol X}) = 1$, then $\Sing {\ol X} = A_1 + A_2$ by \cite{MZ},
but $X$ has two $(-1)$-curves, contradicting Figure 5 in \cite{Ye}.
If $\rho({\ol X}) = 2$, then $\Sing {\ol X} = 2A_1$, and $E, \ell'$ are the
only $(-1)$-curves on $X$ by Figure 6 in \cite{Ye}. So both $(-1)$-curves
are $G$-stable, hence the image on ${\ol X}$ of $E$ is $G$-stable and
contractible, contradicting $\rho({\ol X} \quot G) = 1$.
\par

Suppose $Y \cong \Sigma_0$. We also consider the 
inverse of the morphism $p$. Let $\ell_i$ be the fibres of the two different 
$\BP^1$-fibrations which pass through the fundamental point $P$. After blowing up the 
point $P$, there are two choices of taking the center $Q$ of the next blow-up. 
Namely, if $Q$ is on the proper transform of $\ell_i$ with $i = 1$ say,
then $\ell_2' + \Delta(1) + E + \ell_1'$ has the dual graph in the above
paragraph and we will reach the same contradiction;
if $Q$ does not lie on the proper transforms 
of $\ell_i$, then we obtain Figure $21$ ($E := E_1$). 
We have $\rho({\ol X}) = 3$ and $\Sing {\ol X} = A_1$
(see Figure 5 and 6 in \cite{Ye}). 
We have thus verified all the assertions of Lemma 10.
\QED

We finally consider the case where $\Delta$ has two connected components $\Delta_1$ and 
$\Delta_2$. As in Lemma 2, write
\[
\begin{array}{lll}
{(\Delta_1)}_\red &=&\Delta_1(s)+\cdots+\Delta_1(1) \\
{(\Delta_2)}_\red &=&\Delta_2(t)+\cdots+\Delta_2(1),
\end{array}
\]
where $D.\Delta_1(s)=D.\Delta_2(t)=1$. 

\begin{lem}
Suppose that $D$ is irreducible, that $\Delta_\red$ has two connected components 
$(\Delta_i)_{\rm red}$ ($i = 1,2$) of lengths $s, t$ and that 
$K_X^2 \le 7$. Then the following assertions hold,
where each of Figures $22 \sim 43$
contains the graph of
$f^{-1}{\ol D} + f^{-1}(\Sing {\ol X}) = D + \Delta_{\red} + A$.

\begin{enumerate}
\item[{\rm (0)}] If $s \ne t$ then both $\Delta_i$ are $G$-stable.
\item[{\rm (1)}]
$K_X^2 \ge 3+D^2$, and $K_X^2 \ge 4+D^2$ provided $s=t=1$.
\item[{\rm (2)}]
$2 \le s+t \le 6-D^2$, and $s+t \le 5-D^2$ provided $s=t=1$.
\item[{\rm (3)}]
$-1 \le D^2 \le 4$, and $0 \le D^2 \le 3$ provided $s=t=1$.
\item[{\rm (4)}]
For any $(-1)$-curve $E$ ($\ne D$) on $X$, we have either $E.\Delta_\red
=1$, or $E.\Delta_\red=2, E.\Delta_i(1)=1\ (i=1,2)$ and $-K_X=E+D+\Delta_\red$. One can 
say simply that $D+\Delta_\red+E$ is a simple loop in the latter case.

\item[{\rm (5)}]
In case $D^2=-1$, there are ten possibilities for $D + \Delta$
(see Figures $22 \sim 31$).
In Figure $25$ (resp. $27$), there is an element 
$g$ in $G$ such that $g(E_1) = E_2$ (resp. $g(\Delta_1) = \Delta_2$). 
In Figures $28$ and $31$, no $E_i$ or $F_j$ is $G$-stable.

\item[{\rm (6)}]
In case $D^2=0$, there are nine possibilities (see Figures $32 \sim 40$).
In Figures $34, 36$ and $40$, no $E_i$ or $F_j$ is $G$-stable.
In Figures $33$, $35$ and $39$, there is an element
$g$ in $G$ such that $g(\Delta_1) = \Delta_2$.

\item[{\rm (7)}]
The case $D^2=4$ is impossible. In case $D^2=3$, there is one possibility;
see Figure $41$ where $E$ is the only $(-1)$-curve on $X$.
\item[{\rm (8)}]
The case $D^2=2$ is impossible.
\item[{\rm (9)}]
In case $D^2=1$, there are two possibilities (see Figures $42$ and $43$).
In Figure $42$, the $E_1$ and $E_2$ are the only $(-1)$-curves on $X$.
In Figure $43$), $E, E_1, E_2$ are the only $(-1)$-curves on $X$
and there is an element $g$ in $G$ such that $g(E_1) = E_2$.

\end{enumerate}
\end{lem}

\Proof
The last part of assertion (5) or (6) follows from the arguments for Figure
16 in Lemma 10 and Figure 31 below. Indeed, in Figures 27, 33 and 39,
the argument of Figure 16 shows that $g(E_2) = F_2$ for some $g$ in $G$;
in Figure 35, the argument of Figure 31 shows that
$g(E_1) = F_1$ for some $g$ in $G$.

\svskip \noindent
(0) is clear for $D + \Delta$ is $G$-stable. 

\svskip \noindent
(1)\ \ Since 
\[
K_X^2>2-\Delta.D-D.K_X=4+D^2-\Delta.D=2+D^2+\frac{1}{s+1}+\frac{1}{t+1}
\]
by Lemma 3, it follows that $K_X^2\ge 3+D^2$, and that $K_X^2\ge 4+D^2$ if $s=t=1$. 

\svskip \noindent
(2)\ \ Since 
\[
1+s+t \le s+t+\rho(\ol{X}) \le \rho(X)=10-K_X^2 < 8-D^2-\left(\frac{1}{s+1}+\frac{1}{t+1}\right),
\]
we conclude that $2 \le s+t \le 6-D^2$, and that $s+t \le 5-D^2$ if $s=t=1$.

\svskip \noindent
(3)\ \ By Lemma 8, $D^2 \ge -1$. 
Since $K_X^2 \le 7$, we have $D^2 \le 4$. If $s=t=1$ and $D^2=-1$, then $D+\Delta_\red$ is 
negative semi-definite. Hence $D^2\ge 0$ if $s=t=1$. Furthermore, $D^2 \le 3$ if $s=t=1$ (cf. 
the assertion (1)).

\svskip \noindent
(4)\ \ It follows from Lemmas 3 and 8.

\svskip \noindent
(5)\ \ Let $D^2 = -1$. Then $s+t \ge 3$ by the assertion (3).
We may assume that $s \le t$.  Set $\Gamma_0 := \Delta_1(s)+2D+\Delta_2(t)$.
Then $\Phi := \Phi_{|\Gamma_0|} : X \rightarrow B$ ($B \cong \BP^1$) is a
$\BP^1$-fibration with a singular fibre $\Gamma_0$ and
cross-sections $\Delta_1(s-1)$ (if $s \ge 2$)
and $\Delta_2(t-1)$ (if $t \ge 2$).  Since $\Gamma_0$ is $G$-stable,
$G$ permutes fibres of $\Phi$. We will often use the following:

\svskip \noindent
{\bf Fact 12.}\ \ {\em Suppose $s \ge 2$. There is then a composite of blow downs 
$X \rightarrow \Sigma_2$
of all $(-1)$-curves in fibres not meeting
the cross-section $\Delta_1(s-1)$ so that $\Delta_1(s-1)$ becomes 
the minimal section $M$ and $\Delta_2(t-1)$ becomes a section
disjoint from $M$ and with self-intersection $2$.}
\svskip

Suppose $s = 1$ and $t = 2$.
If $\Gamma_1$ ($\ne \Gamma_0$) is a singular fibre, then it is of type $I_n$ 
with two $(-1)$-curves $E_1, E_2$ as in Lemma 2, and we may assume that the cross-section 
$\Delta_2(1)$ meets $E_2$; then $E_1 \cap \Delta = \emptyset$,
contradicting Lemma 3. So $\Phi$ has only one singular fibre $\Gamma_0$.
Hence $K_X^2=6$. See Figure $22$. 
\par

Suppose $s=1$ and $t \ge 3$. 
Let $\Gamma_1$ be the singular fibre of $\Phi$ containing 
$\Delta_2(t-2)+\cdots+\Delta_2(1)$. By the argument above,
$\Gamma_1$ is the only singular fibre besides $\Gamma_0$.
If $\Gamma_1$ is of type $I_n$ ($n = t-2$) then
no $(-1)$-curve in $\Gamma_1$ is $G$-stable and hence $t = 3$;
see the argument for Figure 16.
Now according to the different types of $\Gamma_1$
in Lemma 2, we have three cases: 
Figure $23$ with $t=3$ and $K_X^2=4$, Figure $24$ with $t=5$ and $K_X^2=3$ and 
Figure $25$ with $t=3$ and $K_X^2=4$. 
\par

Suppose $s=t=2$. According to 
the number of singular fibres and using Fact 12 and assertion (4), 
we have three cases : Figures $26, 27, 28$ all with $K_X^2=3$.
\par

Suppose $s=2$ and $t \ge 3$.  
Let $\Gamma_1$ be the fibre of $\Phi$ containing $\Delta_2(1) + \cdots +
\Delta_2(t-2)$. Then $\Gamma_1 = E_1 + \Delta_2(1) + \cdots + \Delta_2(t-2) + E_2$
is an ordered linear chain and a singular fiber of type $I_{t-2}$ in Lemma 2, 
and we may assume that $E_1$  intersects the cross-section $\Delta_1(1)$.
Note that $\Gamma_1$, $E_i$ are all $G$-stable (see the assertion (0)); so if
$t = 3$ then the image on ${\ol X}$ of $E_2$ is $G$-stable and 
contractible, contradicting $\rho({\ol X}\quot G) = 1$. Thus $t \ge 4$
and $t = 4, 5$. If $\Phi$ has exactly one type $I_0$ singular fibre 
$\Gamma_2 = F_1 + F_2$, then $\Gamma_2$ is $G$-stable;
we may assume $F_1$ (resp. $F_2$) intersects the cross-section 
$\Delta_1(1)$ (resp. $\Delta_2(t-1)$) and hence both $F_i$ are $G$-stable,
but then the image on ${\ol X}$ of $F_1$ is $G$-stable and contractible,
contradicting $\rho({\ol X}\quot G) = 1$.
By the above arguments and by Fact 12 and assertion (4),
we see that $t = 5$, $K_X^2 = 2$ and $\Gamma_0, \Gamma_1$ are the only singular fibres
of $\Phi$. See Figure 29.
\par

Suppose $s=t=3$. We also consider the possible singular fibres of $\Phi$. 
Then by Fact 12 and assertion (4),
two cases given in Figures $30$ and $31$ survive. In Figure 31,
if one of $E_i, F_j$ say $E_1$ is $G$-stable, then 
$2(E_1+\Delta_1(2))+\Delta_1(3)+\Delta_1(1)$ is a $G$-stable
fibre of another $\BP^1$-fibration $\Psi$ where all exceptional divisors
of $f : X \rightarrow {\ol X}$ are contained in fibres; 
this $\Psi$ induces a $\BP^1$-fibration on ${\ol X} \quot G$,
whence $\rho({\ol X}\quot G) \ge 2$, a contradiction.
\par

The remaining case is $s=3$ and $t=4$. In this case we can show by the argument for the
case $s = t = 3$ that there is no possibility.

\svskip \noindent
(6)\ \ Let $D^2=0$. Then $K_X^2\ge 3+D^2 = 3$ and $s+t \le 6-D^2=6$. 
We again assume $s \le t$.  Let $\Phi : X \rightarrow B$ ($B \cong \BP^1$)
be the $\BP^1$-fibration defined by $|D|$, for which 
$\Delta_1(s)$ and $\Delta_2(t)$ are cross-sections. 
Since $D$ is $G$-stable, $G$ permutes fibres of $\Phi$.
\par

Suppose $s = t = 1$. Again, we consider all possibilities of 
the singular fibres of $\Phi$ listed up in Lemma 2. By Fact 12 and the assertion (4),
there are five possibilities : Figures $32 \sim 36$.
\par

Suppose $s=1$ and $t \ge 2$.  Let  $\Gamma_1$  be the singular fibre of $\Phi$  
containing $(\Delta_2)_\red-\Delta_2(t)$. Then 
$\Gamma_1 = E_1 + \Delta_2(1) + \cdots + \Delta_2(t-1) + E_2$ is an
ordered linear chain as in Lemma 2 and we may assume that
$E_1$ meets the cross-section $(\Delta_1)_{\rm red}$.  Note that
$\Gamma_1$, $E_i$ are all $G$-stable.  If $t = 2$, then the image
on ${\ol X}$ of $E_2$ is $G$-stable and contractible, contradicting
$\rho({\ol X}\quot G) = 1$; so $t \ge 3$.
If $t = 3$, then we reach a contradiction as in Figure 31 by
considering another $\BP^1$-fibration $\Phi_1$ defined by
$|2(E_2+\Delta_2(2)) + \Delta_2(3) + \Delta_2(1)|$. Thus $t \ge 4$
and hence $t = 4, 5$. By the argument for the case $D^2= -1, s = 2, t \ge 4$,
it is impossible that $\Phi$ contains a unique type $I_0$ singular fibre.
By the arguments above and by Fact 12 and assertion (4),
we see that $t = 5$, $K_X^2 = 3$ and $\Gamma_1$ is the only singular fibre
of $\Phi$. See Figure 37.
\par

Suppose $s \ge 2$ and $t \ge 2$. Then $(s,t)=(2,2),(2,3), (2,4)$ or $(3,3)$. 
If $s+t=6$, then 
\[
7 \le s+t+\rho(\ol{X}) \le \rho(X)=10-K_X^2 \le 7.
\]
So $\rho(\ol{X})=1$ and $\Sing\ol{X}=A_s+A_t$ with $(s,t)=(2,4),(3,3)$. But these cases are 
impossible by \cite{MZ}. Suppose $(s,t)=(2,3)$. Then we have 
\[
6 \le s+t+\rho(\ol{X}) \le \rho(X) \le 7.
\]
Hence either $\rho(X)=6$ or $\rho(X)=7$. If $\rho(X)=6$, then $\rho(\ol{X})=1$ and $\Sing\ol{X}
=A_2+A_3$ which is impossible by \cite{MZ}. Suppose $\rho(X)=7$. If $\rho(\ol{X})=1$ then 
$\Sing\ol{X}=A_1+A_2+A_3$, which is impossible by \cite{MZ}. If $\rho(\ol{X})=2$ then 
$\Sing\ol{X}=A_2+A_3$, which is impossible by \cite{Ye}. 
So the remaining case is $(s,t)=(2,2)$. 
\par

Let $\Gamma_i$ be the singular fibre of $\Phi$ containing $\Delta_i(1)$.
If $\Gamma_1$ is of type $I_1$ in Lemma 2, then so is $\Gamma_2$, 
and we can write $\Gamma_i = E_i + \Delta_i(1) + F_i$ such that
$F_1$ (resp. $E_2$) meets the cross-section $\Delta_2(2)$
(resp. $\Delta_1(2)$); since $D + \Delta$ is $G$-stable, the image
on ${\ol X}$ of $E_1 + F_2$ is $G$-stable and also contractible,
contradicting $\rho({\ol X} \quot G) = 1$.
The above argument, Fact 12 and the assertion (4) imply that there
are only three possibilities. See Figures $38 \sim 40$, where $K_X^2 = 3$
in all three cases.

\svskip \noindent
(7)\ \ Consider the case $D^2=4$. Then $K_X^2=7$ and $s+t \le 6-D^2=2$. 
Hence $(s,t)=(1,1)$, while then 
$s+t \le 5-D^2=1$, which is absurd. So the case $D^2=4$ does not occur. 
\par

Consider the case $D^2=3$. Then $K_X^2 \ge 3+D^2=6$, whence $K_X^2=6$ or $7$. 
Meanwhile, $s+t \le 6-D^2=3$. So, $(s,t)=(1,1), (1,2)$. 
\par

Suppose $(s,t)=(1,1)$. Then $K_X^2 \ge 4+D^2=7$. Hence $K_X^2=7$ by the 
assumption. Since we have
\[
3=\rho(X) \ge \rho(\ol{X})+s+t \ge 1+2,
\]
we have $\rho(\ol{X})=1$ and $\Sing\ol{X}=2A_1$, which is impossible by \cite{MZ}. 

Suppose $(s,t)=(1,2)$. Then we have 
\[
4 \ge 10-K_X^2 = \rho(X) \ge \rho(\ol{X})+s+t \ge 1+3,
\]
whence follows that $K_X^2=6, \rho(\ol{X})=1$ and $\Sing\ol{X}=A_1+A_2$. 
Note that there is only one $(-1)$-curve $E$ on $X$ and $E . (\Delta_1)_{\red} = 
E . \Delta_2(1) = 1$ by the assertion (4) and by Figure 5 in \cite{Ye}.
See Figure $41$.

\svskip \noindent
(8)\ \ We shall show that $D^2=2$ is impossible. In fact, since $2 \le s+t \le 6-D^2=4$, 
we have $(s,t)=(1,1), (1,2), (1,3)$ or $(2,2)$, where we assume $s \le t$. Note that 
$K_X^2 \ge 3+D^2=5$. If $s+t=4$, we have $\rho(\ol{X})=1$ and $\Sing\ol{X}=A_s+A_t$ with 
$(s,t)=(1,3), (2,2)$. But we cannot find these cases in the table in \cite{MZ}. 
\par

Suppose $(s,t)=(1,1)$. Then $K_X^2 \ge 4+D^2=6$, whence $K_X^2=6,7$. We utilize the inequality 
\[
3=1+s+t \le \rho(\ol{X})+s+t \le \rho(X)=10-K_X^2.
\]
If $K_X^2=7$, we have $\rho(\ol{X})=1$ and $\Sing\ol{X}=2A_1$, which is impossible by 
\cite{MZ}. If $K_X^2=6$, then either $\rho(\ol{X})=1$ and $\Sing\ol{X}=3A_1$, or 
$\rho(\ol{X})=2$ and $\Sing\ol{X}=2A_1$. The former case is impossible by \cite{MZ}. In 
the latter case, take a $(-1)$-curve $E$. If $E.\Delta_\red=1$, say $E.(\Delta_1)_\red=1$, 
the blowing-down of $E+(\Delta_1)_\red$ brings $X$ to $\Sigma_2$, while the image $\wt{D}$ of 
$D$ satisfies $\wt{D}^2=3$, which is impossible on $\Sigma_2$. If $E.\Delta_\red=2$, then 
$E.\Delta=1$, contradicting Lemma 3. 
\par

Suppose $(s,t)=(1,2)$. By an argument similar to the above using the inequalities
\[
4=1+s+t \le \rho(\ol{X})+s+t \le \rho(X)=10-K_X^2 \le 5,
\]
we see that $\Sing\ol{X}=A_1+A_2$ and either $K_X^2=6$ and $\rho(\ol{X})=1$, or
$K_X^2=5$ and $\rho(\ol{X})=2$.  
In the first case, there is only one $(-1)$-curve on $X$ and $E . \Delta_i(1) = 1$
($i = 1, 2$) by the assertion (4) and Figure 5 in \cite{Ye}.
Let $p : X \to Y$ be the blow-down of $E, \Delta_2(1), \Delta_2(2)$.
Since $K_Y^2 = K_X^2 + 3 = 9$, we have $Y \cong {\BP^2}$. However,
$p_*(D)^2 = 3$, which is impossible on ${\BP^2}$.
In the second case, there are only three $(-1)$-curves on $X$
one of which is disjoint from $\Delta$ (see Figure 6 in \cite{Ye}); 
this contradicts Lemma 3.

\svskip \noindent
(9)\ \ Now we treat the case $D^2=1$. Note that $K_X^2 \ge 3+D^2=4$ and $s+t \le 6-D^2 
=5$. We shall prove the following claim.
\svskip

\noindent
{\bf Claim 13.}\ \ {\em There exists a $(-1)$-curve, say $E$, on $X$ such that 
$E.\Delta_\red =2$.}

\svskip \noindent
\Proof
Consider the morphism $q : X \to \BP^2$ defined by $|D|$. Then $D$ is the pull-back of a 
line by $q$. Since $D$ is not touched, $\Delta_1(s)$ and $\Delta_2(t)$ are mapped to lines 
$\ell_1$ and $\ell_2$, respectively. Let $P:=\ell_1\cap \ell_2$. Then $P$ is one of the 
fundamental point of the morphism. We consider to reverse the morphism. Let $E_1$ be the 
$(-1)$-curve appearing by the blowing-up of $P$. If $E_1$ stays as a $(-1)$-curve on $X$, 
then it is a $(-1)$-curve we require for. Otherwise, one of the intersection points 
$P_1, P_2$ of $E_1$ with the proper transforms of $\ell_1$ and $\ell_2$ is blown up, but 
both points are not; if both points are blown up, there will appear a $(-n)$-curve with 
$n \ge 3$, a contradiction. Then the proper transform of $E_1$ on $X$ is contained in 
$\Delta_\red$. Now blow up one of the points $P_1, P_2$ and apply the same argument as above 
to the $(-1)$-curve $E_2$ appearing from the blow-up. We have just only to continue 
this argument. 
\QED

In case $s+t=5$, we have $6 \le s+t+\rho(\ol{X}) \le \rho(X)=10-K_X^2 \le 6$. Hence 
$\rho(\ol{X})=1$ and $\Sing\ol{X}=A_s+A_t$ with $(s,t)=(1,4), (2,3)$. But these cases are 
impossible by \cite{MZ}. 
\par

Suppose $s+t=4$. Then we have 
\[
5 \le s+t+\rho(\ol{X}) \le \rho(X)=10-K_X^2 \le 6.
\]
If $K_X^2=5$, then $\rho(\ol{X})=1$ and $\Sing\ol{X}=A_s+A_t$ with $(s,t)=(1,3), (2,2)$. 
These cases do not exist by \cite{MZ}. If $K_X^2=4$, by \cite{MZ}, \cite{Ye}, 
we have $(s, t) = (1, 3)$, and either $\rho(\ol{X})=1$ and 
$\Sing\ol{X}=2A_1+A_3$ 
or $\rho(\ol{X})=2$ and $\Sing\ol{X}=A_1+A_3$.
See Figures $42$ and $43$, where $E_1$ or $E$ is as in Claim 13.
The part about the uniqueness of the $(-1)$-curves in the
assertion (9) follows from Figures 5 and 6 in \cite{Ye}.
Since $G$ acts on the set of $(-1)$-curves on $X$, in
Figure 43, $E$ is $G$-stable and each element of $G$ either
stabilizes or switches $E_1$ and $E_2$; so the existence of
$g$ in $G$ with $g(E_1) = E_2$ follows from $\rho({\ol X} \quot G) = 1$
(see the argument for Figure 16).
\par

Suppose $(s,t)=(1,1)$. Then a $(-1)$-curve $E$ as in Claim 13 
has $E.\Delta=1$, contradicting Lemma 3. 
\par

Suppose $(s,t)=(1,2)$. Consider the $\BP^1$-fibration $\Phi : X \to B$ defined 
by $|\Gamma_0|$ where $\Gamma_0 := \Delta_1(1)+2E+\Delta_2(1)$.
Then we can make the Hirzebruch surface $\Sigma_2$ out of $X$ with the image of 
$\Delta_2(2)$ as the minimal section. The blow-down of $E, \Delta_1(1)$ increases 
$D^2$ by $1$. Since $D.\Delta_2(2)=1$, the blow-down of $E, \Delta_1(1)$ is not enough 
to bring $X$ to $\Sigma_2$. Hence there exists a singular fibre $\Gamma_1$ of
type $I_n$ which then contains a $(-1)$-curve 
$E_1$ meeting the cross-section $D$. 
This is a contradiction by Lemma 3. This ends the proof of Lemma 11.

\section{Determination of the group $G$ action on $X$}

In this section, we shall consider all $43$ triplets $({\ol X}, {\ol D}; G)$ 
in Theorem A in the introduction, determine the action of the finite group 
$G$ on $X$ and give examples.

\svskip \noindent
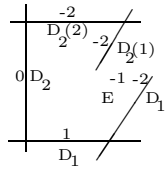
{\bf Fig 1}. Let $\psi : X \to {\BP^2}$ be the blow-down of 
$E+\Delta_2(1)+\Delta_2(2)$ to a point $P$. 
Clearly, there is an induced faithful action of $G$ on $\BP^2$ 
such that $\psi$ is $G$-equivariant. So $G$ is a subgroup
of $\PGL_2(\C)$ stabilizing each component of the triangle
$\psi(\Delta_1+D_1+D_2)$. We may assume that the three vertices
of the triangle are at $[1, 0, 0], [0, 1, 0], [0, 0, 1]$. 
Then $G$ is a subgroup of $\{\diag[1, b, c] \mid bc \ne 0\} 
\subseteq \PGL_2(\C)$. Conversely, 
any finite subgroup of $\{\diag[1, b, c] \mid bc \ne 0\}$
can act faithfully on this ${\ol X}$ fitting Figure 1, 
such that $\rho({\ol X} \quot G) = 1$ 
(noting that $\rho({\ol X}) = 1$ already).

\par \svskip \noindent
{\bf Fig 2 $\sim$ 6}.  Let $H$ be the (normal) subgroup of $G$ stabilizing $D_1$ 
(and hence also $D_2$), and let $g$ be an element in $G$ switching $D_1$ and $D_2$
(see Lemma 7). Then $G = \langle g, H \rangle$.
Note that $H$ is abelian. This is because
at the point $D_1 \cap D_2$, all elements of $H$ can be diagonalized
simultaneously with the same eigenvectors along the directions of $D_1$ and $D_2$.

\par
In Figure 2, one can show that $H$ is cyclic. Indeed, $H$ fixes the three intersection
points of the cross-section $\Delta_1(2)$ with the three singular fibres of
different types, and hence $H|_{\Delta_1(2)} = \{$id $\}$.
So every $h$ in $H$ is diagonalized as $(1, c_h)$
at $D_1 \cap \Delta_1(2)$ with the common eigenvectors
along the directions of $\Delta_1(2)$ and $D_1$.
Thus $H$ can be embedded in $\C^*$ via $h \mapsto c_h$ and is cyclic.

\par
Since $\rho({\ol X} \quot \langle g \rangle) = 1$
can be easily checked, we have always $\rho({\ol X} \quot G) = 1$
so long $G$ exists. Note that $G = \langle g \rangle \cong {\Z}/(2)$
is realizable in all these 5 cases (Lemma 7).

\par \svskip \noindent
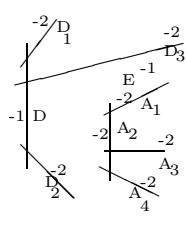
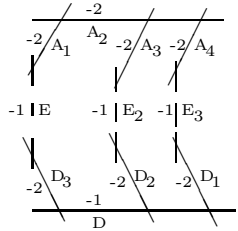
{\bf Fig 7}. Since $\rho({\ol X}) = 1$ and Sing ${\ol X} = 3A_1 + D_4$, 
there are exactly three $(-1)$-curves $E, E_2, E_3$ on $X$ fitting Figure 7', where
there is a ${\BP^1}$-fibration $\Psi$ on $X$ such that
$\Gamma_0' := 2E+A_1+(\Delta_3)_{\red}, 
\Gamma_1' := 2E_2 + A_3 + (\Delta_2)_{\red}, 
\Gamma_2' := 2E_3 + A_4 + (\Delta_1)_{\red}$ 
are all the singular fibres and $A_2$ and $D$ are cross-sections of $\Psi$. 
Clearly, $G$ permutes fibres of
$\Psi$. Let $\psi : X \to {\BP^2}$ be the blow-down of 
$E+A_1, E_2+A_3, E_3+A_4, D$ to 4 points $P_1, \dots, P_4$, respectively. 
Then $\psi$ is $G$-equivariant; $G$ fixes $P_4$ and permutes 
$P_1, P_2, P_3$. We may assume that $\psi(A_2) = \{Z = 0\}$ which is $G$-stable, and 
$P_4 = [0, 0, 1]$ which is $G$-fixed. Let $H$ be the (normal) subgroup of 
$G$ fixing all three points $P_1, P_2, P_3$ 
(and hence fixes the line $\{Z = 0\}$), then $H = \langle h_1 \rangle$ 
for some $h_1 = \diag[1, 1, c_1]$ of order $n_1$.
\par 

Let $\iota : G \rightarrow \Aut\{P_1, P_2, P_3\} = S_3$ be the natural
homomorphism. Then $\im(\iota) = S_3$, ${\Bbb Z}/(3)$, ${\Bbb Z}/(2)$ or $(1)$.
If an element $h_3$ in $G$ acts transitively on the set 
$\{P_1, P_2, P_3\}$, then $h_3^3$ acts trivially on the line $\{Z = 0\}$
and hence $h_3 = \diag[1, \omega, c_3]$, where $\omega = 
\exp(2 \pi \sqrt{-1}/3)$, after the normalization that $h_3$ fixes two points
$[1, 0, 0], [0, 1, 0]$ on the line $\{Z = 0\}$.
If there is further an element $h_2$ in $G$ acting as an involution on
the set $\{P_1, P_2, P_3\}$, one may assume that $h_2(P_1) = P_1$
and $P_1 = [1, 1, 0]$. One can show that 
$h_2 = \begin{pmatrix} 0 & 1 & 0 \\
          1 & 0 & 0 \\
          0 & 0 & c_2 \end{pmatrix}$,
by using the following conditions:
$h_2(P_i) = P_i$ ($i = 1, 4$), $h_2(P_2) = P_3, h_2(P_3) = P_2$,
$P_2 = h_3^j(P_1) = [1, \omega^j, 0], P_3 = h_3^{2j}(P_1)$
($j = 1$ or $2$).
\par 

Suppose that $\im(\iota) = \Z/(2)$ and let $h_2'$ be in $G$ acting as   
an involution on the set $\{P_1, P_2, P_3\}$. Then $h_2' = \diag[1, -1, c_2]$
after the normalization that $h_2'$ fixes two points $[1, 0, 0], [0, 1, 0]$ 
on the line $\{Z = 0\}$.
\par 

Replacing $h_3$ (resp. $h_2$ or $h_2'$) by its power we may assume that
ord$(h_3) = 3^{n_3}$ (resp. ord$(h_2)$ or ord$(h_2')$ is $2^{n_2}$).
Thus either $G = \langle h_1, h_2, h_3 \mid h_1h_i = h_ih_1 (i=2,3), h_2^2, h_3^3, 
h_2^{-1}h_3h_2h_3 \in \langle h_1 \rangle \rangle$, or 
$G = \langle h_1, h_3 \mid h_1h_3 = h_3h_1, h_3^3 \in \langle h_1 \rangle \rangle$, or
$G = \langle h_1, h_2' \mid h_1h_2' = h_2'h_1, (h_2')^2 \in \langle h_1 \rangle \rangle$, or
$G =\langle h_1 \rangle$. We have an exact sequence:
$$
(1) \to \langle h_1 \rangle \to G \to G/\langle h_1 \rangle \to (1)
$$
where $G/\langle h_1 \rangle = S_3$, $\Z/(3)$, $\Z/(2)$ or $(1)$.
We may take $G = (1)$ since then $\rho({\ol X} \quot G) = \rho({\ol X}) = 1$.

\svskip \noindent
In {\bf Figures 8 $\sim$ 13} below, let $Q_i := D \cap (\Delta_i)_{\red}$. Let
$q_1 : X_1 \to X$ be the blow-up of a point $R_2$ ($\ne Q_i$) on $D$ 
with $J_0$ the exceptional curve (we may choose $R_2$ to be a $G$-fixed 
point if it exists, but always $R_2 \ne Q_i$). Then $-K_{X_1} = J_0 + q_1'(2D+
(\Delta_1)_{\red} + (\Delta_2)_{\red} + (\Delta_3)_{\red})$, 
which is nef and big. By the Riemann-Roch theorem and 
the Kawamata-Viehweg vanishing theorem, dim $|-K_{X_1}| = 1$.
Let $q_0 : X_0 \to X_1$ be the blow-up of the unique base point of 
$|-K_{X_1}|$ (which must lie on $J_0$)
(cf. Proposition 2 at page 40 of \cite{De}, 
or Lemma 1.7 in \cite{DZ}) with $O$ the exceptional divisor.
Then there is an elliptic fibration $\gamma : X_0 \to {\BP^1}$ with $O$ 
the zero section and $T_0 := 2D+J_0+(\Delta_1)_{\red}+(\Delta_2)_{\red}+
(\Delta_3)_{\red}$ as a singular fibre (we use, by the abuse of notations, the same symbol like 
$\Delta_i$ to denote its proper transform on $X_0$) which is of type 
$I_0^*$. Let $\Aut_0(X_0) = \{g \in \Aut(X_0) \mid g(O) = O\}$. Then there is an
induced $\Aut_0(X_0)$ action on $X$ so that
$q = q_1 \circ q_0 : X_0 \to X$ is $\Aut_0(X_0)$-equivariant.
Clearly, $\Aut_0(X_0) = \{g \in \Aut(X) \mid g(R_2) = R_2\}$.
Let $T$ be a general fibre of the elliptic fibration on $X_0$. 
Then $\Aut_0(X_0)$ contains $\Aut_0(T) = \{g \in \Aut(T) \mid g$ fixes the point
$O \cap T \}$ $\cong \Z/(m)$ with $m = 2, 4$ or $6$ (see Cor. 4.7 in \cite{Ha} at page 321).
Hence $\Aut_0(T)$ (and $\Aut_0(X_0)$) contains an involution $\sigma : t \mapsto -t$.
\par 

Let $\iota : G \to \Aut\{Q_1, Q_2, Q_3\} = S_3$ be the natural homomorphism. 
Let $H := \Ker(\iota)$, which then acts trivially on $D$; 
in particular, $H \subseteq \Aut_0(X_0)$. At the point $Q_1$, every element $h$ in $H$ 
has the directions of $D$ and $\Delta_1$ as eigenvectors with respect to the 
eigenvalues $1, \lambda_h$. So $H$ can be embedded into $k^*$ and hence $H$ is cyclic. 
Note that we have an exact sequence
$$
(1) \to H \to G \to G/H \to (1)
$$
where $G/H = (1)$, $\Z/(2)$, $\Z/(3)$ or $S_3$.
Let $G_0$ be any finite group of $\Aut_0(X_0)$ containing the involution 
$\sigma$ of a general fibre; note that $\sigma$ is in the centre of 
$Aut_0(X_0)$. We shall show that in each of Figures 8 $\sim$ 13, we can take $G = G_0$
so that $\rho({\ol X} \quot G) = 1$.

\svskip \noindent
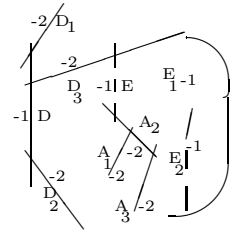
{\bf Fig 8}. In this case the elliptic fibration $\gamma$ has $T_0$, 
$T_1 = B + A_1 + A_2 + A_3$ which is of type $I_4$, and a few irreducible fibres 
as singular fibers. As in Figure 10 below, considering the height pairing,
we can show that $\sigma(E) = E$, $\sigma(E_1) = E_2$ and
$\rho({\ol X} \quot G_0) = 1$.

\par \svskip \noindent
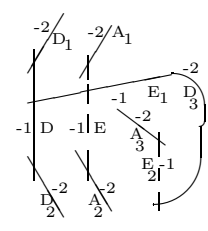
{\bf Fig 9}. In this case, $\gamma$ has $T_0$, $T_i = A_i + B_i$ ($i = 1, 2, 3$) 
each of which is of type $I_2$ or $III$, and a few irreducible fibres as singular fibers.
As in Figure 10 below, considering the height pairing,
we can show that $\sigma(E) = E$, $\sigma(E_1) = E_2$ and
$\rho({\ol X} \quot G_0) = 1$.

\par \svskip \noindent
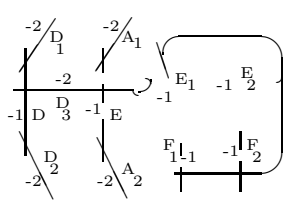
{\bf Fig 10}. In this case, $\gamma$ has $T_0$, $T_i = A_i + B_i$ ($i = 1,2$),
each of which is of type $I_2$ or $III$, and a few irreducible fibres as singular fibers.
On the surface $X_0$, by the height pairing in \cite{Sh}, 
$\langle E, E \rangle = 2 \chi({\SO}_{X_0}) + 2 E . O - (1 + 1/2 + 1/2) =
0$, whence $E$ is a torsion in ${\rm MW}(\gamma)$;
one can see easily that $E$ is a $2$-torsion and hence
$\sigma(E) = E$;
also $\langle E_i, E_i \rangle = 1 = \langle F_i, F_i \rangle$, 
$\langle E_i, F_i \rangle = \chi({\SO}_{X_0}) + E_i . O + F_i . O - 
E_i . F_i - (1 + 0 + 0) = -1$, $\langle E_i, E_j \rangle = 0 =
\langle E_i, F_j \rangle$ for $i \ne j$.
On the surface $X$, since $\sigma$ stabilizes the fibre 
$2E + A_1 + A_2$ of $\Phi$, it permutes fibres of $\Phi$, whence
$\sigma(E_i) = E_j$ or $F_j$ for some $j$.
Note that in $MW(\gamma)$, $E_i + \sigma(E_i) = 0$
and hence $\langle E_i + \sigma(E_i), E_i + \sigma(E_i) \rangle = 0$.
By the calculation above, we must have $\sigma(E_i) = F_i$.
On the surface $X$, since $\Pic X$ is generated over 
$\Q$ by the fibre components and a 2-section $\Delta_3$, the $\Pic {\ol X}$
is generated over $\Q$ by the images ${\ol E}$,
${\ol E}_i$, ${\ol F}_j$ of $E, E_i, F_j$ with $2{\ol E} =
{\ol E}_i + {\ol F}_i$. So it follows that $\rho({\ol X} \quot G_0) = 1$.

\par \svskip \noindent
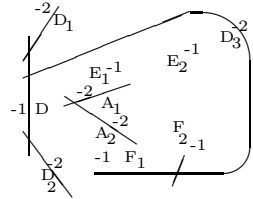
{\bf Fig 11}. In this case, $\gamma$ has $T_0$, $T_1 = A_1 + A_2 + B$,
which is of type $I_3$ or $IV$, and a few irreducible fibres as singular fibers.
On the surface $X_0$, we can calculate the height pairing
and find that $\langle E_i + F_i, E_i + F_i \rangle = 0$;
so $E_i+F_i$ is torsion and it must be zero in $MW(\gamma)$
for the latter is torsion free by \cite{Pe} or \cite{Mi}.
So $\sigma(E_i) = -E_i = F_i$ in $MW(\gamma)$.
As in Figure 10, $\rho({\ol X} \quot G_0) = 1$.

\par \svskip \noindent
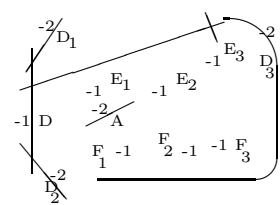
{\bf Fig 12}. In this case, $\gamma$ has $T_0$, $T_1 = A + B$
which is of type $I_2$ or $III$, and a few irreducible fibres as singular fibers.
As in Figure 11, $\sigma(E_i) = F_i$ ($1 \le i \le 3$) 
and hence as in Figure 10, $\rho({\ol X} \quot G_0) = 1$.

\par \svskip \noindent
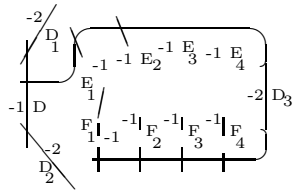
{\bf Fig 13}. In this case, $\gamma$ has $T_0$ and a few irreducible fibres 
as singular fibers. As in Figures 10 and 11, $\sigma(E_i) = F_i$ ($1 \le i \le 4$) 
and $\rho({\ol X} \quot G_0) = 1$.

\par \svskip \noindent
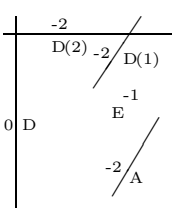
{\bf Fig 14}. Let $\psi : X \to {\BP^2}$ be the blow-down of 
$E + \Delta(1) + \Delta(2)$ to a point $P$. We may assume that
$\psi(A) = \{X = 0\}$ and $\psi(D) = \{Y = 0\}$ so that 
$P = [0, 0, 1]$. Then $\psi$ is $G$-equivariant and
$G \subseteq \{g = (a_{ij}) \in \PGL_2(\C) \mid a_{21} = 0\ \mbox{and $g$
is lower triangular}\}$. Conversely, any finite group in $\PGL_2(\C)$ of such form
can act on this ${\ol X}$ faithfully so that $\rho({\ol X} \quot G) = 1$.
(Note that $\rho({\ol X})$ is already $1$.)

\par \svskip \noindent
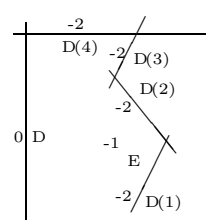
{\bf Fig 15}. By blowing down $E$, we are reduced to Figure 14.
Let $P' = E \cap \Delta(2)$ which is infinitely near to the point $P$
defined in Figure 14. Then $G \subseteq \{g = (a_{ij})\in \PGL_2(\C) \mid g(P') = P',
a_{21} = 0\ \mbox{and $g$ is lower triangular}\}$. Conversely, any finite group in 
$\PGL_2(\C)$ of such form can act on this ${\ol X}$ faithfully so that 
$\rho({\ol X} \quot G) = 1$. (Note that $\rho({\ol X})$ is already $1$.)

\par \svskip \noindent
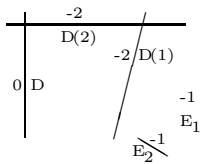
{\bf Fig 16}. Let $g$ be the element in $G$ switching $E_1$ and $E_2$
(see Lemma 10), and let $H$ be the (normal) subgroup of $G$ stabilizing $E_1$
(and hence $E_2$). Let $\psi : X \to {\BP^2}$ be the blow-down
of $E_2 + \Delta(1) + \Delta(2)$, which is $H$-equivariant.
As in Figure 14, $H \subseteq \{h = (a_{ij})\in \PGL_2(\C) \mid a_{21} = 0\ 
\mbox{and $h$ is lower triangular}\}$. Note that $H$ is normal in $G$ 
and $G = \langle g, H \rangle$.
\par 

Here is an example where $G = \langle g \rangle$ and ord$(g) = 2$.
Let $\Sigma_4$ be the Hirzebruch surface with the $(-4)$-curve $M$
and a section $B$ disjoint from $M$. Let
$Y \to \Sigma_4$ be the blow-up of a point not on $M$.
Note that $M + B$ is 2-divisible in the Picard lattice.
Let $X \rightarrow Y$ be the double cover branched along
$M+B$ with $\langle g \rangle = \Gal(X/Y)$. Let $\Delta(2)$ be
the inverse on $X$ of $M$, $\Delta(1)$ the proper transform on $X$ 
of the fibre through the centre of the blow-up $Y \rightarrow \Sigma_4$,
$E_1+E_2$ the inverse of the exceptional curve of the same blow-up,
and let $D$ be the inverse of any general fibre. Then 
Figure 16 appears on this $X$ so that $\rho({\ol X} \quot \langle g \rangle) = 1$,
where $X \to {\ol X}$ is the blow-down of $\Delta(1)+\Delta(2)$.

\par \svskip \noindent
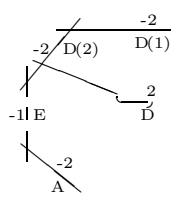
{\bf Fig 17}. Since $E$ is the only $(-1)$-curve on $X$ 
(see Figure 5 in \cite{Ye}), $E$ is $G$-stable.
Let $\psi : X \to {\BP}^2$ be the blow-down of
$E+\Delta(2)+\Delta(1)$ to a point say $P := [0, 0, 1]$. Then
$\psi$ is $G$-equivariant. We may also assume that $\psi(A) = \{X = 0\}$.
Note that $\psi(D)$ is a conic touching $A$ at $P$ with order $2$.
Then $G \subseteq \{g = (a_{ij})\in \PGL_2(\C) \mid g(\psi(D)) = \psi(D), 
\ \mbox{$g$ is lower triangular}\}$. Conversely, 
any finite group in $\PGL_2(\C)$ of such form
can act on this ${\ol X}$ faithfully so that $\rho({\ol X} \quot G) = 1$.
(Note that $\rho({\ol X})$ is already $1$.)

\par \svskip \noindent
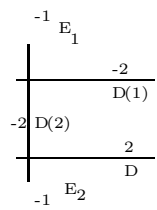
{\bf Fig 18}. Note that $E_1, E_2$ are the only $(-1)$-curves on $X$
(see Figure 6 in \cite{Ye}) and hence $G$ stabilizes $E_1+E_2$.
Let $g$ be the element in $G$ switching $E_1$ and $E_2$ (see Lemma 10)
and let $H$ be the (normal) subgroup of $G$ stabilizing $E_1$ (and hence $E_2$).
Let $\psi : X \to {\BP^2}$ be the blow-down of $E_2 + \Delta(2) + \Delta(1)$
to a point say $P := [0, 0, 1]$. Then $\psi$ is $H$-equivariant.
$H$ stabilizes the line $\psi(E_1)$ defined by $\{X = 0\}$ say, and also
the conic $\psi(D)$ touching $\psi(E_1)$ at $P$ with order $2$.
As in Figure 17, $H \subseteq \{h = (a_{ij})\in \PGL_2(\C) | h(\psi(D)) = \psi(D),
\ \mbox{$h$ is lower triangular}\}$. Note that $H$ is normal in $G$
and $G = \langle g, H \rangle$.
\par 

Here is an example with $G = \langle g \rangle$ and $\ord(g) = 2$.
Let $X \to Y$ and $\Delta(i), E_i$ 
be as in Figure 16, but we let $D$ be the inverse
on $X$ of $B$.

\par \svskip \noindent
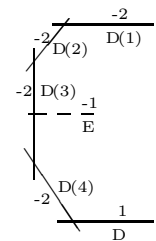
{\bf Fig 19}. Since $E$ is the unique $(-1)$-curve on $X$
(see Figure 5 in \cite{Ye}), it is $G$-stable.
By blowing down $E$, we are reduced to Figure 20.
Set $P' := E \cap \Delta(3)$ which is an infinitely near point
of the point $P$ in Figure 20.
Thus $G \subseteq \{g = (a_{ij})\in \PGL_2(\C) \mid g(P') = P',
\ \mbox{$a_{ij} \ne 0$ only when $i = j$ or $(i,j) = (3,1)$} \}$.
Conversely, any finite group in $\PGL_2(\C)$ of such form
can act on this ${\ol X}$ faithfully so that $\rho({\ol X} \quot G) = 1$.
(Note that $\rho({\ol X})$ is already $1$.)

\par \svskip \noindent
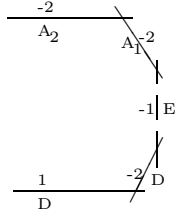
{\bf Fig 20}. Since $E$ is the only $(-1)$-curve on $X$ (see Figure 5 in
\cite{Ye}), it is $G$-stable. Let $\psi : X \to {\BP^2}$ be the blow-down
of $E+A_1+A_2$ to a point, say $P := [0, 1, 0]$. Then $\psi$ is $G$-equivariant.
$G$ fixes $P$ and $[0, 0, 1]$ which is the intersection of
two $G$-stable lines $\psi(\Delta) = \{X = 0\}$ and $\psi(D) = \{Y = 0\}$ say, 
whence $G \subseteq \{(a_{ij})\in \PGL_2(\C) \mid a_{ij} \ne 0\ \mbox{only when $i = j$ 
or $(i,j) = (3,1)$} \}$.
Conversely, any finite group in $\PGL_2(\C)$ of such form
can act on this ${\ol X}$ faithfully so that $\rho({\ol X} \quot G) = 1$.
(Note that $\rho({\ol X})$ is already $1$.)

\par \svskip \noindent
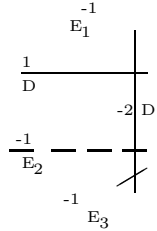
{\bf Fig 21}. Let $\psi : X \rightarrow {\BP^2}$ be the blow-down of
the $E_i$ to the points $P_i$. We can show that $E_i$ are the only
$(-1)$-curves on $X$ and hence $\psi$ is $G$-equivariant.
So $G \subseteq \Aut_S = \{g \in \PGL_2(\C) \mid g(S) = S, g({\hat D}) = {\hat D} \}$, 
where $S = \{P_1, P_2, P_3\}$ is a subset on the line $\psi(\Delta)$
and ${\hat D}$ ($= \psi(D)$) is a second line.
Since $\rho({\ol X} \quot G) = 1$,
the $G$ acts on $S$ transitively. Conversely,
any finite group in $\Aut_S$ acting transitively on $S$
can act on this ${\ol X}$ faithfully so that $\rho({\ol X} \quot G) = 1$.

\par \svskip \noindent
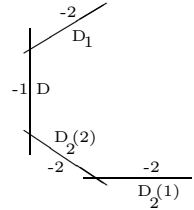
{\bf Fig 22}. Let $\psi : X \to {\BP^2}$ be the blow-down of 
$D + \Delta_2(2) + \Delta_2(1)$ to a point say $P := [0,0,1]$. Then
$\psi$ is $G$-equivariant. $G$ fixes $P$ and stabilizes the line
$\psi(\Delta(1)) = \{X = 0\}$ say. Then
$G \subseteq \{g = (a_{ij})\in \PGL_2(\C) \mid \ \mbox{$g$ is lower triangular}\}$.
Conversely, any finite group in $\PGL_2(\C)$ of such form
can act on this ${\ol X}$ faithfully so that $\rho({\ol X} \quot G) = 1$.
(Note that $\rho({\ol X})$ is already $1$.)

\par \svskip \noindent
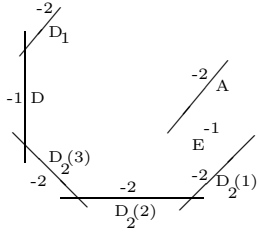
{\bf Fig 23}. Let $\psi : X \rightarrow {\BP^2}$ be the blow-down of
$E+A$ and $D + \Delta_2(3) + \Delta_2(2)$ to points $P := [0, 1, 0]$
and $[0, 0, 1]$ say. Then $\psi$ is $G$-equivariant.
$G$ fixes $P$ and stabilizes two lines $\psi(\Delta_2(1) = \{X = 0\}$
and $\psi(\Delta_1)) = \{Y = 0\}$ say. Then
$G \subseteq \{(a_{ij})\in \PGL_2(\C) \mid \ \mbox{
$a_{ij} \ne 0$ only when $i = j$ or $(i,j) = (3,1)$}\}$.
Conversely, any finite group in $\PGL_2(\C)$ of such form
can act on this ${\ol X}$ faithfully so that $\rho({\ol X} \quot G) = 1$.
(Note that $\rho({\ol X})$ is already $1$.)

\par \svskip \noindent
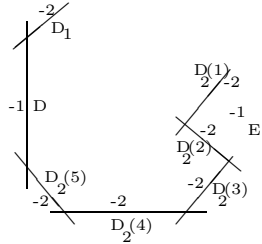
{\bf Fig 24}. By blowing down $E$, we are reduced to Figure 23.
Set $P' := E \cap \Delta_2(2)$ which is an infinitely near point of
the point $P$ in Figure 23. Then 
$G \subseteq \{g = (a_{ij})\in \PGL_2(\C) \mid g(P') = P', \ \mbox{
$a_{ij} \ne 0$ only when $i = j$ or $(i,j) = (3,1)$} \}$.
Conversely, any finite group in $\PGL_2(\C)$ of such form
can act on this ${\ol X}$ faithfully so that $\rho({\ol X} \quot G) = 1$.
(Note that $\rho({\ol X})$ is already $1$.)

\par \svskip \noindent
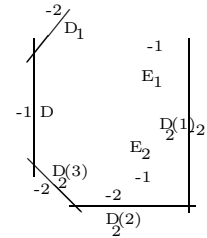
{\bf Fig 25}. Let $g$ be an element in $G$ switching $E_1$ and $E_2$
(see Lemma 11) and let $H$ be the (normal) subgroup of $G$ stabilizing $E_1$
(and hence $E_2$). Let $\psi : X \rightarrow {\BP^2}$ be the blow-down
of $D + \Delta_1$ and $E_2+\Delta_2(1)+\Delta_2(2)$ to two points
$P_1 = [0, 1, 0]$ and $P_2 = [0, 0, 1]$ say. Then $\psi$ is $H$-equivariant.
$H$ fixes $P_i$ and stabilizes the two lines
$\psi(\Delta_2(3)) = \{X = 0\}$ and $\psi(E_1) = \{Y = 0\}$
say. Then $H \subseteq \{(a_{ij})\in \PGL_2(\C) \mid 
a_{ij} \ne 0\ \mbox{only when $i = j$ or $(i,j) = (3,1)$} \}$.
Note that $H$ is normal in $G$ and $G = \langle g, H \rangle$.
\par 

Here is an example where $G = \langle g \rangle$ and ord$(g) = 2$.
Let $M$ be the $(-2)$-curve on the Hirzebruch surface $\Sigma_2$,
$B$ a section disjoint from $M$ and $L_i$ are two distinct fibres.
Let $p : Y \to \Sigma_2$ be the blow-up of a point on $L_2$ other than
$L_2 \cap M$, the point $L_1 \cap M$ and its infinitely near point
lying on the proper transform of $M$;
let ${\hat E}$, ${\hat D}$, ${\hat \Delta}_1$ be irreducible curves
on $Y$ which are (the proper transforms of)
the corresponding exceptional curves.
Since $M+B$ is 2-divisible in the Picard lattice, there is a
double cover $X \to Y$ branched along ${\hat D} + p'(M+B)$
with $\langle g \rangle = \Gal(X/Y)$.
Let $D$, $\Delta_1$, $\Delta_2(1)$, $\Delta_2(2)$, $\Delta_2(3)$, $E_1+E_2$
be the strict inverses of $\wh{D}$, $L_1$, $L_2$, $\wh{\Delta}_1$ and $\wh{E}$, 
respectively. Then Figure 25 appears on this $X$ so that 
$\rho({\ol X} \quot \langle g \rangle) = 1$,
where $X \to {\ol X}$ is the blow-down of $\Delta_1+ \sum_j \Delta_2(j)$.

\par \svskip \noindent
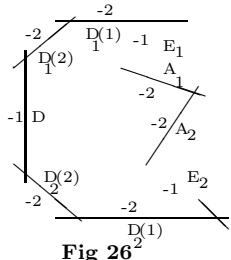
{\bf Fig 26}. Let $H$ be the (normal) subgroup of $G$ stabilizing $\Delta_1$
(and hence also all of $\Delta_i(j)$).
Blowing down $D$, the Figure becomes Figure 5, whence
$H$ is abelian. Thus either $G = H$, or $G = \langle g, H \rangle$
where  $g$  switches $\Delta_1$ and $\Delta_2$.

\par 
Let $\psi : X \to {\BP^2}$ be the blow-down of $D + \Delta_2(2) + \Delta_2(1)$
and $E_1 + A_1 + A_2$ to $P_1 = [1, 0, 0]$ and $P_2 = [0, 1, 0]$ say.
Then $\psi$ is $H$-equivariant. $H$ fixes $P_3 = \Delta_1(1) \cap \Delta_1(2)$
with coordinates, say $[0, 0, 1]$, on ${\BP^2}$ and also two points
$P_i$ ($i = 1, 2$) on the line $\psi(E_2) = \{Z = 0\}$. 
Thus $H \subseteq \{\diag[1, b, c] \mid bc \ne 0 \}$.

\par
Conversely, each finite group $G = H$ in $\PGL_2(\C)$ of the form above
or $G = \langle g \rangle \cong \Z/(2)$ is realizable as a group of 
automorphisms on this ${\ol X}$ such that $\rho({\ol X} \quot G) = 1$. 
(Indeed, $\rho({\ol X}) = 1$ already; see Lemma 7 for the second case.)

\par \svskip \noindent
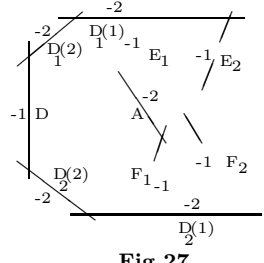
{\bf Fig 27}. Let $g$ be in $G$ switching $\Delta_1$ and $\Delta_2$ 
(Lemma 11). Let $H$ be the (normal) subgroup of $G$ stabilizing $\Delta_1$ 
(and hence $\Delta_2$). As in Figure 2, we have 
$H = \langle h \rangle$ and $G = \langle g, h \rangle$ with
$h|_{\Delta_i(1)} =$id. The case $G = \langle g \rangle \cong \Z/(2)$
is realizable (see Lemma 7; indeed Figure 27 is different from Figure 2 only
in labelling).

\par \svskip \noindent
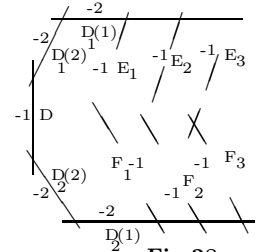
{\bf Fig 28}. Let $H$ be the (normal) subgroup of $G$ stabilizing $\Delta_1$ 
(and hence also $\Delta_2$). Let $H_1$ be the (normal) subgroup of $G$ stabilizing $E_1$ 
and $E_2$ (and hence all $E_i, F_j, \Delta_i(j)$). As in Figure 2, $H$ is abelian and
$H_1 = \langle h_1 \rangle$ with $h_1|_{\Delta_i(1)} =$id.
Note that $G/H \le \Z/(2)$ and $|H/H_1| \le 3$; indeed, $H/H_1$ is abelian and acts on
the set $\{E_1, E_2, E_3\}$. By Lemma 11, either $G/H = \Z/(2)$
or $G = H$ and $H/H_1 = \Z/(3)$.
Each of the case $G = G/H = \Z/(2)$ and the case $G = H$ with $H/H_1 = \Z/(3)$ 
is realizable as a group of automorphisms on this ${\ol X}$ 
such that $\rho({\ol X} \quot G) = 1$ (see Figure 7 and Lemma 7, 
noting that the Figure becomes Figure 4 after the blow down of $D$).

\par \svskip \noindent
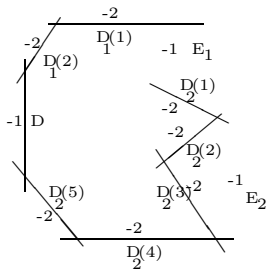
{\bf Fig 29}. Let $\psi : X \to {\BP^2}$ be the blow-down of
of $D + \Delta_2(5) + \Delta_2(4)$, $E_1 + \Delta_2(1) + \Delta_2(2)$
and $E_2$ to points say $[1, 0, 0], [0, 1, 0]$ and $[1, 1, 0]$
on the same line $\psi(\Delta_2(3)) = \{Z = 0\}$.
Then $\psi$ is $G$-equivariant. $G$ fixes these three points and
also the intersection point $\Delta_1(1) \cap \Delta_1(2)$ with coordinates
say $[0, 0, 1]$ on ${\BP^2}$. Then $G = \langle g \rangle$
where $g = \diag [1, 1, c]$.
Conversely, any finite cyclic group can act on this ${\ol X}$ 
faithfully so that $\rho({\ol X} \quot G) = 1$. 
(Note that $\rho({\ol X})$ is already $1$.)

\par \svskip \noindent
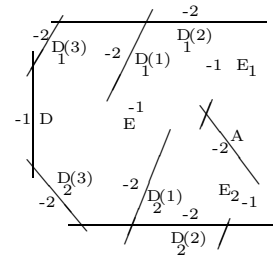
{\bf Fig 30}. Blowing down $D$ which is $G$-stable, the Figure becomes Figure 2.
So either $G = \langle g, h \rangle$ or $G = \langle h \rangle$,
where $h|_{\Delta_i(2)} =$id and $g$ switches $\Delta_1$ and $\Delta_2$.
Each of $G = \langle g \rangle \cong \Z/(2)$ and $G = \langle h \rangle$
is realizable as a group of automorphisms on this ${\ol X}$ 
such that $\rho({\ol X} \quot G) = 1$. 
(Note that $\rho({\ol X}) = 1$ already.)

\par \svskip \noindent
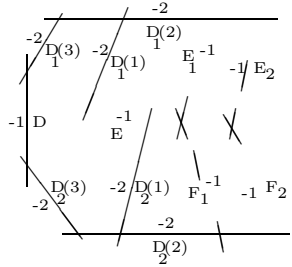
{\bf Fig 31}. Let $H$ be the (normal) subgroup of $G$ stabilizing $\Delta_1$ 
(and hence also $\Delta_2$). Let $H_1$ be the (normal) subgroup of $G$ stabilizing $E_1$ 
(and hence all $E_i, F_j, \Delta_i(j)$). As in Figure 2, $H$ is abelian and
$H_1 = \langle h_1 \rangle$ with $h_1|_{\Delta_i(2)} =$id.
As in Figure 28, by Lemma 11, either $G/H = \Z/(2)$,
or $G = H$ and $H/H_1 = \Z/(2)$.
Each of the case $G = G/H \cong \Z/(2)$ and the case $G = H$ with $H/H_1 = \Z/(2)$ 
is realizable as a group of automorphisms on this ${\ol X}$ 
such that $\rho({\ol X} \quot G) = 1$ (see Lemma 7 and Figure 7).

\par \svskip \noindent
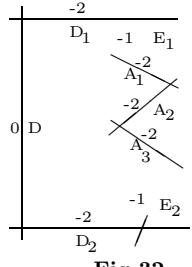
{\bf Fig 32}.  Let $\psi : X \rightarrow {\BP^1} \times {\BP^1}$
be the blow-down of $E_1+A_1$ and $E_2+A_3$ to two points $P_1$ and $P_2$,
respectively. Then $\psi$ is $G$-equivariant.
Thus $G$ is a subgroup of $\Aut_S(\BP^1\times\BP^1) := \{g \in \Aut(\BP^1 \times \BP^1) \mid
g(S) = S, g(\wh{D}) = \wh{D} \}$, where $S = \{P_1, P_2\}$ is a set of two points on the same 
fibre of the first ruling and $\wh{D}$ ($= \psi(D)$)
is a second fibre of the same ruling.
Let $\iota : G \to \Aut(S) = \Z/(2)$ be the natural homomorphism
and let $H = \Ker(\iota)$. Then all elements of $H$ can be diagonalized
simultaneously at $P_1$ with the same eigenvectors along the directions of
$\psi(\Delta_1)$ and $\psi(A_2)$. Thus $H \subseteq \{\diag[b, c] \mid bc \ne 0 \}$.
Note that $G/H = (0)$ or ${\Bbb Z}/(2)$.
Conversely, any finite subgroup $G$ of $Aut_S(\BP^1\times\BP^1)$ can act on
${\ol X}$ faithfully with $\rho({\ol X} \quot G) = 1$.
(Note that $\rho({\ol X}) = 1$ already.)

\par \svskip \noindent
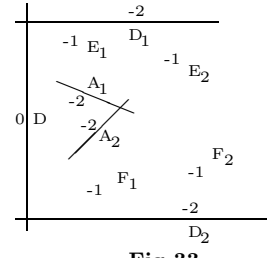
{\bf Fig 33}. Let $g$ be in $G$ switching $\Delta_1$ and $\Delta_2$ (cf. Lemma 11).
As in Figure 2, $G = \langle g, h \rangle$, where $h|_{\Delta_i} =$id,
and the case $G = \langle g \rangle \cong \Z/(2)$
is realizable (see Lemma 7; indeed Figure 33 is different from Figure 5 only
in labelling).

\par \svskip \noindent
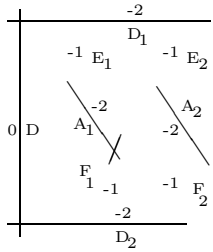
{\bf Fig 34}. Let $\psi : X \rightarrow {\BP^1} \times {\BP^1}$
be the blow-down of $E_1, E_2, F_1, F_2$ to four points $e_1, e_2, f_1, f_2$,
respectively. Then $\psi$ is $G$-equivariant.
Let $\iota : G \to Aut\{e_1, e_2, f_1, f_2\} = S_4$ be the natural
homomorphism. Then $\im(\iota)$ is contained in the Klein
fourgroup $V = \langle (e_1e_2)(f_1f_2), (e_1f_1)(e_2f_2) \rangle$ of $S_4$.
We assert that $\im(\iota) \subseteq \langle (e_1f_2)(e_2f_1) \rangle$
is impossible.
Indeed, if this assertion is false, then $G$ permutes fibres of the
${\BP^1}$-fibration $\Psi$ with singular fibres $2E_1 + A_1 + \Delta_1$,
$2F_2 + A_2 + \Delta_2$, where all components of $f^{-1}(\Sing {\ol X})$
are contained in fibres of $\Psi$; this leads to $\rho({\ol X} \quot G) \ge 2$
as in Lemma 11, which is a contradiction. So the assertion is true.
Note that $G$ is a subgroup of $Aut_S(\BP^1\times\BP^1) := \{g \in \Aut(\BP^1\times \BP^1) \mid
g(S) = S, g(\wh{D}) = \wh{D} \}$, where $S = \{e_1, e_2, f_1, f_2\}$ is the intersection
of four fibres, two from each ruling and $\wh{D}$ ($= \psi(D)$)
is a fifth fibre. By the assertion, we have:

\par \svskip \noindent
$(*)$ $G|_S$ equals either the Klein group $V$ or $\langle (e_1e_2)(f_1f_2) \rangle$
or $\langle (e_1f_1)(e_2f_2) \rangle$.

\par \svskip \noindent
Let $\iota : G \to \Aut\{D \cap \Delta_1, D \cap \Delta_2\} = \Z/(2)$
be the natural homomorphism. As in Figure 32, we have
$H := \Ker(\iota) \subseteq \{\diag[b, c] \mid bc \ne 0 \}$,
and $G/H = (0)$ or $\Z/(2)$.
Conversely, any finite subgroup $G$ of $Aut_S(\BP^1\times\BP^1)$ satisfying $(*)$ can act on
${\ol X}$ faithfully with $\rho({\ol X} \quot G) = 1$.

\par \svskip \noindent
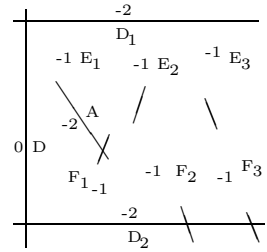
{\bf Fig 35}. Let $g$ be in $G$ switching $\Delta_1$ and $\Delta_2$ (Lemma 11).
Let $H$ be the (normal) subgroup of $G$ stabilizing $\Delta_1$ 
(and hence also $\Delta_2$). As in Figure 6, $H$ is abelian,
$G = \langle g, H \rangle$, and the case $G = \langle g \rangle \cong \Z/(2)$
is realizable (see Lemma 7; indeed Figure 35 is different from Figure 6 only
in labelling).

\par \svskip \noindent
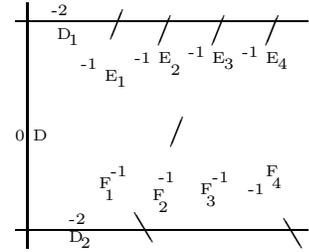
{\bf Fig 36}. Let $H_1$ (resp. $H_2$) be the (normal) subgroup of $G$
stabilizing all $E_i$ (resp. stabilizing $\Delta_1$). Then both $H_i$
are normal in $G$ such that $G/H_2 = (0)$ or ${\Bbb Z}/(2)$
and $H_2/H_1 \subseteq S_4$. As in Figure 2, $H_1 = \langle h \rangle$.
If $G = H_2$, then $\rho({\ol X} \quot G) = 1$
implies that $H_2$ acts on the set $\{E_1, \dots, E_4\}$
transitively, i.e., $H_2/H_1$ is a transitive subgroup of $S_4$.
Conversely, $\langle g \rangle \cong {\Bbb Z}/(2)$ can actually act
on ${\ol X}$ such that $\rho({\ol X} \quot \langle g \rangle) = 1$
and $g(E_i) = F_i$ for all $i$ (see Lemma 7, noting that Figure 36 is
different from Figure 4 only in labelling).

\par \svskip \noindent
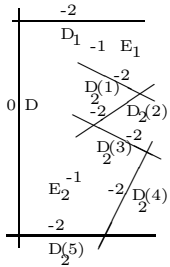
{\bf Fig 37}. Let $\psi : X \to {\BP^2}$ be the blow-down of 
$E_2 + \Delta_2(4) + \Delta_2(5)$ and $E_1 + \Delta_2(1) + \Delta_2(2)$
to points, say $P_1 = [0, 0, 1]$ and $P_2 = [0, 1, 0]$. 
Then $\psi$ is $G$-equivariant. $G$ fixes the $P_i$ and the intersection
of $D$ and $\Delta_1$ with coordinates say $[1, 0, 0]$ on ${\BP^2}$.
Set $P_1' := E_2 \cap \Delta_2(4)$. Then $G \subseteq
\{g = \diag[1, b, c] \in \PGL_2(\C) \mid g(P_1') = P_1' \}$.
Conversely, any finite group in $\PGL_2(\C)$ of such form
can act on this ${\ol X}$ faithfully so that $\rho({\ol X} \quot G) = 1$.
(Note that $\rho({\ol X})$ is already $1$.)

\par \svskip \noindent
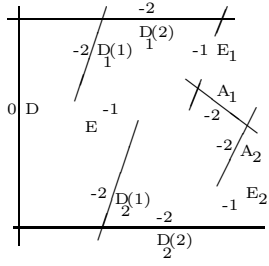
{\bf Fig 38}. Let $H$ be the (normal) subgroup of $G$ stabilizing $\Delta_1$ 
(and hence also $\Delta_2$). As in Figure 2, $H = \langle h \rangle$,
where $h|_{\Delta_i(2)} =$id. Note that $G/H \le \Z/(2)$.
The case $G = G/H \cong \Z/(2)$ is realizable (see Lemma 7; 
indeed, blowing down $E$, the Figure becomes Figure 5).

\par \svskip \noindent
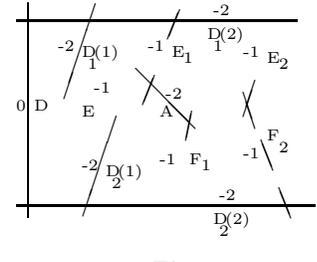
{\bf Fig 39}. Let $g$ be in $G$ switching $\Delta_1$ and $\Delta_2$ (Lemma 11).
As in Figure 2, $G = \langle g, h \rangle$, where $h|_{\Delta_i(2)} =$id,
and the case $G = \langle g \rangle \cong \Z/(2)$
is realizable (see Lemma 7; indeed Figure 39 is different from Figure 2 only
in labelling).

\par \svskip \noindent
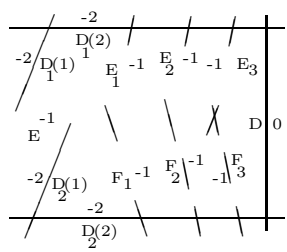
{\bf Fig 40}. Figure 40 is identical with Figure 28 with only difference in labelling.

\par \svskip \noindent
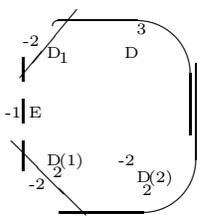
{\bf Fig 41}. Let $\psi : X \to {\BP^2}$ be the blow-down of
$E+\Delta_2(1)+\Delta_2(2)$ to a point $P_1$. Then $\psi$ is
$G$-equivariant. $G$ fixes $P_1$ and the intersection point
$P_2$ of $D$ with $\Delta_1$.
Thus $G \subseteq \{g \in \PGL_2(\C) \mid
g(P_i) = P_i (i=1,2), g(\wh{D}) = \wh{D} \}$,
where $\wh{D}$ ($= \psi(D)$) is a conic intersecting
the line $L_{P_1P_2}$ ($= \psi(\Delta_1)$) at the points $P_i$.
Conversely, any finite group in $\PGL_2(\C)$ of such form
can act on this ${\ol X}$ faithfully so that $\rho({\ol X} \quot G) = 1$.
(Note that $\rho({\ol X})$ is already $1$.)

\par \svskip \noindent
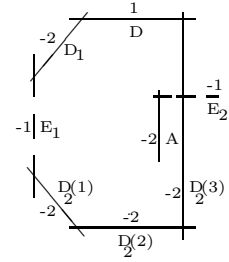
{\bf Fig 42}. Let $\psi : X \to {\BP^2}$ be the blow-down of
$E_1+\Delta_2(1)+\Delta_2(2)$ and $E_2+A$ to points
$P_2 = [1, 0, 0]$ and $P_3 = [1, 1, 0]$ say. Then 
$\psi$ is $G$-equivariant. $G$ fixes three points
$P_i$ ($i = 1,2,3$) where $P_1 = D \cap \Delta_2(3)$ 
with $P_1 = [0, 1, 0]$ say, all lying on the line 
$\psi(\Delta_2(3)) = \{Z = 0\}$, and also the intersection point
$P_4 = D \cap \Delta_1$ with $P_4 = [0, 0, 1]$ say.
Then $G = \langle g \rangle$ with $g = [1, 1, c]$.
Conversely, any finite cyclic group
can act on this ${\ol X}$ faithfully so that $\rho({\ol X} \quot G) = 1$.
(Note that $\rho({\ol X})$ is already $1$.)

\par \svskip \noindent
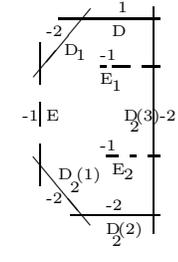
{\bf Fig 43}. Let $\psi : X \to {\BP^2}$ be the blow-down of
$E_1$, $E_2$ and $E + \Delta_2(1) + \Delta_2(2)$ to points
$P_i$ ($i = 1,2,3$), respectively. Then $\psi$ is $G$-equivariant.
$G$ fixes three points $P_3 = [1, 0, 0]$ say, $P_4 := D \cap \Delta_2(3)
= [0, 1, 0]$ and $P_5 := D \cap \Delta_1 = [0, 0, 1]$. 
Note that $P_1, \dots, P_4$ lie on the same $G$-stable 
line $\psi(\Delta_2(3)) = \{Z = 0\}$ say.
So $G \subseteq \{\diag[1, b, c] \mid bc \ne 0\}$. 
Let $g$ be an element in $G$ switching $E_1$ and $E_2$ (see Lemma 11).
Then $g$ switches $P_1$ and $P_2$. 
We may assume that
$P_1 = [1, 1, 0]$. Now $g(P_1) = P_2$ and $g(P_2) = P_1$ imply
that $P_2 = [1, -1, 0]$ and $g = [1, -1, c_1]$.
Let $H$ be the (normal) subgroup of $G$ fixing $P_1$ (and hence $P_2$).
Then $H = \langle h \rangle$ for some $h = \diag[1, 1, c_2]$.
Thus $G = \langle g = \diag[1, -1, c_1], h = \diag[1, 1, c_2] \rangle$.
Conversely, any finite group in $\PGL_2(\C)$ of such form
can act on this ${\ol X}$ faithfully so that $\rho({\ol X} \quot G) = 1$.

\par \svskip \noindent
Theorem A is a consequence of the lemmas in \S 1.
Theorem B is proved in the arguments above. For instance, the assertion that
$\kappa({\ol X} \setminus {\ol D}) = -\infty$ in the Hypothesis (H),
follows from the observation that 
$-(K_X + D + \Delta) = -f^*(K_{\ol X} + {\ol D})$ is nef and big.
Indeed, from the construction of the action of $G$ on ${\ol X}$,
we see that $\rho({\ol X} \quot G) = 1$. So the $G$-stable divisor
$-(K_{\ol X} + {\ol D})$ is either numerically trivial, or ample or anti-ample
(see Lemmas 1 and 2).
Now the observation that $-(K_{\ol X} + {\ol D}) . {\ol D} =
-(K_X + D + \Delta) . D = 2 - D . \Delta > 0$ shows that
$-(K_{\ol X} + {\ol D})$ is ample and hence $-(K_X + D + \Delta)$
is nef and big.

\par \svskip \noindent
Theorem C follows from the classification of the group $G$ in \S 2.
Indeed, for $K_X^2 \le 4$, we see that either $G$ is a subgroup
of $\PGL_2$ as in Theorem C, or there is a sequence of subgroups
of $G$ such that the factor groups are abelian, or $G$ is as in
the case of Figure 25. The easy calculation of $K_X^2$
is given below.

\par \vskip 1pc \noindent
{\bf Lemma 14}
For the $X$ in Figure $m$, we calulate $K_X^2$.
\begin{enumerate}
\item[{\rm (1)}]
$K_X^2 = 2$, if $m$ is one of $7 - 13$, $29 - 31$.
\item[{\rm (2)}]
$K_X^2 = 3$, if $m$ is one of $2 - 3$, $24$, $26 - 28$, $37 - 40$.
\item[{\rm (3)}]
$K_X^2 = 4$, if $m$ is one of $4 - 6$, $23$, $25$, $32 - 36$, $42-43$.
\item[{\rm (4)}]
$K_X^2 = 5$, if $m$ is one of $15$, $19$.
\item[{\rm (5)}]
$K_X^2 = 6$, if $m$ is one of $1$, $14$, $16 - 18$, $20 - 22$, $41$.
\end{enumerate}

\par \vskip 2pc \noindent
Masayoshi Miyanishi

\par \noindent
Department of Mathematics
\par \noindent
Graduate School of Science
\par \noindent
Osaka University 
\par \noindent
Toyonaka, Osaka 560-0043
\par \noindent
Japan
\par \noindent
e-mail : miyanisi$@$math.sci.osaka-u.ac.jp

\par \vskip 1pc \noindent
De-Qi Zhang
\par \noindent
Department of Mathematics
\par \noindent
National University of Singapore
\par \noindent
2 Science Drive 2
\par \noindent
Singapore 117543
\par \noindent
Republic of Singapore
\par \noindent
e-mail : matzdq$@$math.nus.edu.sg


\par \vskip 2pc

\setlength{\unitlength}{1973sp}%
\begingroup\makeatletter\ifx\SetFigFont\undefined%
\gdef\SetFigFont#1#2#3#4#5{%
  \reset@font\fontsize{#1}{#2pt}%
  \fontfamily{#3}\fontseries{#4}\fontshape{#5}%
  \selectfont}%
\fi\endgroup%



\begin{thebibliography}{25}
\bibitem{Bo}
E. Bombieri, Canonical models of surfaces of general type, 
Inst. Hautes �udes Sci. Publ. Math. Vol. {\bf 42} (1973), 171--219.

\bibitem{Br}
E. Brieskorn, Rationale Singularit\"{a}ten komplexer Fl\"{a}chen, Invent. Math. {\bf 4}
(1967/1968), 336--358.

\bibitem{De}
M. Demazure, S\'eminaire sur les singularit\'es des surfaces, Lecture Notes in Mathematics, 
Vol. {\bf 777}, Springer, Berlin, 1980.

\bibitem{DZ}
I. Dolgachev and D. -Q. Zhang, Coble rational surfaces, Amer. J. Math. {\bf 123} (2001), 
79--114. 

\bibitem{Ha}
R. Hartshorne, Algebraic Geometry, Graduate Texts in Mathematics, Vol. {\bf 52}
(1977). 

\bibitem{Ka}
Y. Kawamata, Crepant blowing-up of 3-dimensional canonical singularities and 
its application to degenerations of surfaces, 
Ann. of Math. {\bf 127} (1988), 93--163. 

\bibitem{KMM}
Y. Kawamata, K. Matsuda, and K. Matsuki, Introduction to the minimal model
problem, Algebraic geometry, Sendai, 1985, 283--360, 
Adv. Stud. Pure Math. Vol. {\bf 10}, North-Holland, Amsterdam-New York, 1987. 

\bibitem{KZ}
J. Keum and D. -Q. Zhang, Algebraic surfaces with quotient singularities -
including some discussion on automorphisms and fundamental groups, 
Proc. Alg. Geom. in East Asia, A. Ohbuchi et al. (eds), 3-10 Aug. 2001, Japan.

\bibitem{Ko}
J. Kollar et al., Flip and abundance for algebraic threefolds, Asterisque
{\bf 211} (1992).

\bibitem{Ko2}
J. Kollar, Shafarevich maps and plurigenera of algebraic varieties,
Invent. Math. {\bf 113} (1993), 177--215.

\bibitem{Mi}
R. Miranda, Persson's list of singular fibers for a rational elliptic surface,
Math. Z. {\bf 205} (1990), 191--211.

\bibitem{M1}
M. Miyanishi, Non-complete algebraic surfaces, Lecture Notes in Mathematics, 
Vol. {\bf 857}, Springer-Verlag 1981.

\bibitem{oas}
M. Miyanishi, Open algebraic surfaces, CRM Monograph Series, {\bf 12}, 
Amer. Math. Soc.

\bibitem{KM}
M. Miyanishi and K. Masuda, Open algebraic surfaces with finite group actions, 
Transform. Groups {\bf 7} (2002), 185--207.

\bibitem{MZ}
M. Miyanishi and D.-Q. Zhang, Gorenstein del Pezzo surfaces, I; II; J. of Alg.
{\bf 118} (1988), 63--84; {\bf 156} (1993), 183--193.

\bibitem{Og}
K. Oguiso, Families of hyperkahler manifolds, math.AG/{\bf 9911105}.
    
\bibitem{OS}
K. Oguiso and T. Shioda, The Mordell-Weil lattice of a rational 
elliptic surface, Comment. Math. Univ. St. Paul. {\bf 40} (1991), 
83--99.

\bibitem{Pe}
U. Persson, Configurations of Kodaira fibers on rational elliptic surfaces,
Math. Z. {\bf 205} (1990), 1--47. 

\bibitem{Sh}
T. Shioda, On the Mordell-Weil lattices, 
Comment. Math. Univ. St. Paul. {\bf 39} (1990), 211--240.

\bibitem{Ye}
Q. Ye, On Gorenstein log del Pezzo surfaces, Japan. J. Math. 28 (2002), 87 -- 136.

\bibitem{Z1}
D. -Q. Zhang, Logarithmic del Pezzo surfaces with rational double
and triple singular points, Tohoku Math. J. {\bf 41} (1989), 399--452.

\bibitem{Z2}
D. -Q. Zhang, Automorphisms of finite order on rational surfaces 
(with an Appendix by I. Dolgachev), J. of Alg. {\bf 238} (2001), 560--589. 

\bibitem{Z3}
D. -Q. Zhang, Automorphisms of finite order on Gorenstein del Pezzo surfaces,
Trans. Amer. Math. Soc. {\bf 354} (2002), 4831--4845.

\end{thebibliography}
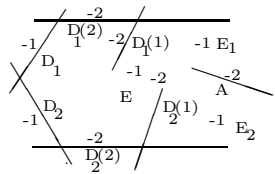
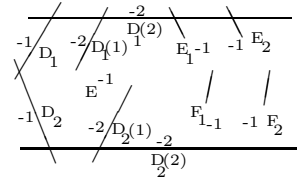
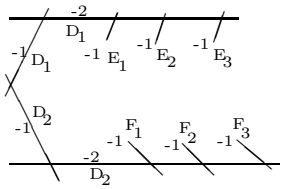
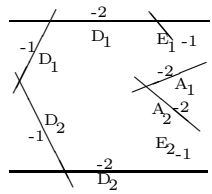
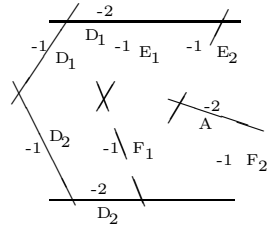
\end{document}